\documentclass[reqno,12pt]{amsart}
\usepackage{amsmath,amsthm,amssymb,amsfonts,amscd}
\usepackage{graphicx,color}
\usepackage{boxedminipage}

\input xy
\xyoption{all}

\newcommand{\R}{{\mathbb{R}}}
\newcommand{\Z}{{\mathbb{Z}}}
\newcommand{\C}{{\mathbb{C}}}

\newcommand{\T}{{\mathbb{T}}}

\newcommand{\ba}{\begin{array}}
\newcommand{\ea}{\end{array}}

\newcommand{\bp}{\begin{pmatrix}}
\newcommand{\ep}{\end{pmatrix}}

\newcommand{\bps}{\begin{smallmatrix}}
\newcommand{\eps}{\end{smallmatrix}}

\newcommand{\ti}{\tilde}

\newcommand{\la}{\langle}
\newcommand{\ra}{\rangle}

\def \ii{{\bf i}}

\def \cA{{\mathcal A}}

\def \cC{{\mathcal C}}

\def \cF{{\mathcal F}}

\def \dim{\mathrm{dim}}
\def \Ext{\mathrm{Ext}}

\def \Ob{\mathrm{Ob}}

\def \Tri{\mathit{Tr}}
\def \Tw{\mathit{Tw}}

\def \ov#1{\frac{1}{#1}}

\def \delm#1{{\delta_{[#1]}}}

\def \fpartial#1{\frac{\partial}{\partial {#1}}}

\def \wt#1{{\widetilde {#1}}}

\def \f{{\frak f}}

\def \l{{\frak l}}
\def \m{{\frak m}}

\def \l({\left(}
\def \r){\right)}

\def \0{{\bf 0}}
\def \1{{\bf 1}}

\def \nat{\mathrm{tan}^{-1}}
\def \Lag{\frak{Lag}}

\newtheorem{thm}{Theorem}[section]

\newtheorem{cor}[thm]{Corollary}
\newtheorem{lem}[thm]{Lemma}
\newtheorem{prop}[thm]{Proposition}

\theoremstyle{definition}

\newtheorem{defn}[thm]{Definition}

\newtheorem{rem}[thm]{Remark}

\begin{document}

\title{Fukaya categories of two-tori revisited}

\author{Hiroshige Kajiura}
\address{Faculty of Science, Chiba university, 
1-33 Yayoi-cho, Inage-ku, Chiba 263-8522, Japan}
\email{kajiura@math.s.chiba-u.ac.jp}

\begin{abstract}
We construct an 
$A_\infty$-structure of the Fukaya category explicitly 
for any flat symplectic two-torus.
The structure constants of the non-transversal $A_\infty$-products
are obtained as derivatives of those of transversal $A_\infty$-products. 
\end{abstract}

\footnote{This work 
is supported by Grant-in Aid for Scientific Research (C) (18K03293) 
of the Japan Society for the Promotion of Science. }

\renewcommand{\thefootnote}{\arabic{footnote}}
\setcounter{footnote}{0}

 \maketitle

{\small



}

\tableofcontents

 \section{Introduction}
 \label{sec:intro}
 
In this paper, we construct an 
$A_\infty$-structure of the Fukaya category explicitly 
for any flat symplectic two-torus. 
For the objects of the Fukaya category, 
we consider pairs of affine Lagrangian submanifolds 
and $U(1)$ local systems on them. 
On a sequence of objects where adjacent Lagrangians 
intersect transversally, the $A_\infty$-product 
is defined explicitly in \cite{PoZa, Poli:higher}. 
It is defined by counting polygons surrounded by those Lagrangians 
following the original definition in \cite{fukayaAinfty} 
with a modification by \cite{kont:hms}. 
Such transversal $A_\infty$-products are 
discussed also for flat symplectic higher dimensional tori \cite{fukaya:abelian}. 
The purpose of this paper is to present the whole $A_\infty$-products 
explicitly including non-transversal $A_\infty$-products. 
Our strategy is to derive the $A_\infty$-structure 
by applying homological perturbation theory to 
the DG category of holomorphic vector bundles on the mirror dual 
complex tori. 
Thus, our result at the same time implies homological mirror symmetry 
\cite{kont:hms} 
of two-tori. 
This strategy is based on the idea in \cite{KoSo:torus} 
on homological mirror symmetry for torus fibrations or 
a more explicit formulation in \cite{hk:hpt-hms, hk:fukayadeform} which 
are closer to our construction. 
On the other hand, in \cite{hk:Ainftyplane} we give an explicit $A_\infty$-structure 
on the Fukaya category of lines in $\R^2$ by a similar strategy. 
As mentioned in \cite{hk:Ainftyplane}, 
since $\R^2$ is the $\Z^2$ covering space of a torus, 
an explicit $A_\infty$-structure on a flat symplectic torus is 
obtained from the one on $\R^2$ in \cite{hk:Ainftyplane} 
by a few modifications. 
In this paper, we present the result explicitly 
and mention later how to obtain the result from the one
in \cite{hk:Ainftyplane}.
An interesting point is that 
the structure constants of the non-transversal $A_\infty$-products
are obtained as derivatives of those of transversal $A_\infty$-products
(Lemma \ref{lem:multi-beta} (iii)). 

The existence of this explicit $A_\infty$-structure is 
already announced in \cite{hk:rims07nc}. 
A full subcategory of such a Fukaya $A_\infty$-category
consisting of finitely many objects is presented in \cite{hk:hpt-hms}. 
We had planed to present the full result
as we cited [38] in the reference \cite{hk:rims07nc},
but have postponed it for a while. 
Now, we again feel that
it is worth presenting the explicit form. 
Though the result is quite complicated, 
it is useful since the 
$A_\infty$-structure is minimal.
For instance, 
it can be applied to visualizing an isomorphism between
a Lagrangian surgery of two intersecting Lagrangians and
the cone (extension) of a morphism associated to the intersecting point(s) 
of the given Lagrangians. 
Via the homological mirror symmetry, it turns out that 
the cone of the morphism between two holomorphic vector bundles
associated to given two Lagrangians are isomorphic to 
the holomorphic vector bundle associated to 
a Lagrangian surgery of the given Lagrangians. 
Such a cone is recently discussed 
in \cite{kobayashi:2tori,chan-suen17}. In particular, 
in \cite{kobayashi:2tori}, 
the isomorphism is constructed explicitly for two holomorphic vector bundles
$E_1,E_2$ such that $\dim\Ext(E_1,E_2)=1$. 
There, quite complicated calculations show that it actually forms
an isomorphism. 
However, since our Fukaya $A_\infty$-category is minimal, 
we can easily construct such an isomorphism as we demonstrate it 
in subsection \ref{ssec:m3}. 

The relation of cones and Lagrangian surgeries as above
leads to show the homological mirror symmetry of two tori. 
In \cite{abouzaid-smith10},
it is shown as an equivalence of split-closed triangulated categories. 
Our result instead turns to show the homological mirror symmetry
at the level of $A_\infty$-enhancements of such triangulated categories.

This paper is organized as follows. 
In section \ref{sec:Ainfty},
we recall cyclicity in an $A_\infty$-category and 
the construction of triangulated category from
a unital $A_\infty$-category in order to fix our conventions. 
In section \ref{sec:Fukaya},
we present explicitly
the minimal cyclic $A_\infty$-structure for the Fukaya category of a torus. 
It is constructed so that it is $A_\infty$-quasi-isomorphic to
the DG category of the holomorphic vector bundles, i.e.,
the homological mirror symmetry holds. 
In section \ref{sec:hms}, we explain the outline of this construction,
which consists of
a few modification of the construction of the Fukaya $A_\infty$-structure
for symplectic planes \cite{hk:Ainftyplane}. 
In section \ref{sec:auto}, we regard the structure constants
of the $A_\infty$-products as generalizations of theta functions and
discuss a few properties they have. 
Finally, in section \ref{sec:exmp}, 
we discuss a few examples of these structure constants. 
In subsection \ref{ssec:m2}, we calculate the structure constants of
all transversal products $m_2$. 
In subsection \ref{ssec:m3},
we calculate certain non-transversal triple products which are related to
an exact triangle in $\Tri(\cC)$. 
This corresponds to the mirror dual of the exact triangle
discussed in \cite{kobayashi:2tori}. 

Throughout this paper, we denote $\sqrt{-1}=:\ii$.


 \section{Cyclic $A_\infty$-categories}
\label{sec:Ainfty}

In order to fix the conventions including the signs, 
in subsection \ref{ssec:cyclicAinfty}
we recall cyclic $A_\infty$-categories
and related notions. For details, see for instance \cite{hk:rims07nc}. 
Later we construct explicitly the Fukaya category on a two-torus
as a (minimal) cyclic $A_\infty$-category. 

We also include in subsection \ref{ssec:triAinfty}
a quick review of constructing
a triangulated category $\Tri(\cC)$
from a strictly unital $A_\infty$-category $\cC$ 
proposed by \cite{kont:hms}. 
We employ these terminologies only
in subsection \ref{ssec:SL2Z} and \ref{ssec:m3},
so the reader can skip this subsection once. 

In this section, any vector space is that over a fixed field $K$. 
Later we employ the tools recalled in this section with $K=\C$.

 \subsection{$A_\infty$-categories and cyclicity}
\label{ssec:cyclicAinfty}
 
An $A_\infty$-category is a natural extension of an 
$A_\infty$-algebra introduced by J.~Stasheff \cite{jds:hahI,jds:hahII}. 
\begin{defn}[$A_\infty$-category \cite{fukayaAinfty}]
An $A_\infty$-category $\cC$ 
consists of a set of objects $\Ob(\cC)$, a graded vector space
$\cC(a,b)=\oplus_{r\in\Z} \cC^r(a,b)$ 
for each objects $a,b\in\Ob(\cC)$ and a collection of multilinear maps 
\begin{equation*}
  \m:=\{m_n:\cC(a_1,a_2)\otimes\cdots\otimes \cC(a_n,a_{n+1})
  \to \cC(a_1,a_{n+1})\}_{n\ge 1}
\end{equation*}
of degree $(2-n)$ satisfying 
\begin{equation}
 \begin{split}
&0=\sum_{k+l=n+1}\sum_{j=0}^{k-1}
(-1)^\star \\
&\ m_k(w_{12},\cdots,w_{j(j+1)},m_l(w_{(j+1)(j+2)},\cdots,w_{(j+l)(j+l+1)}),
 w_{(j+l+1)(j+l+2)},\cdots,w_{n(n+1)})\ ,
 \end{split}
 \label{Ainfty2}
\end{equation}
for homogeneous elements $w_{i(i+1)}\in\cC(a_i,a_{i+1})$, 
where $\star=(j+1)(l+1)+l(|w_{12}|+\cdots+|w_{j(j+1)}|)$. 
 \label{defn:Ainfty2}
\end{defn}
\begin{defn}[Strict units]
For an $A_\infty$-category $\cC$, an element $1_a\in\cC^0(a,a)$ 
is called the strict unit if it satisfies
\begin{equation*}
 m_2(1_a,w_{ab})=w_{ab},\qquad m_2(w_{ba},1_a)=w_{ba}
\end{equation*}
for any $b\in\Ob(\cC)$, $w_{ab}\in\cC(a,b)$, $w_{ba}\in\cC(b,a)$, and
\begin{equation*}
 m_k(\dots, 1_a,\dots) =0
\end{equation*}
for $k=1$ and $k=3,4,5,\dots$. 
\end{defn}
\begin{defn}[Cyclic $A_\infty$-category]
For an  $A_\infty$-category $\cC$, 
a cyclic structure $\eta$ of degree $|\eta|$ 
is a nondegenerate graded symmetric bilinear map 
\[
\eta:\cC(a,b)\otimes\cC(b,a)\to K
\] 
of degree $-|\eta|\in\Z$ for any $a,b\in\Ob(\cC)$ 
satisfying 
\begin{equation}\label{cyclic-cd}
 \begin{split}
 & \eta(m_n(w_{12},\cdots,w_{n(n+1)}),w_{(n+1)1})
 =(-1)^\star 
 \eta(m_n(w_{23},\cdots,w_{(n+1)1}),w_{12})\ , \\
 & \star = {n+(|w_{23}|+\cdots +|w_{(n+1)1}|)|w_{12}|}
 \end{split}
\end{equation}
for each $n\ge 1$. Then, $(\cC,\eta)$, ot simply $\cC$,
is called a cyclic $A_\infty$-category 
of degree $|\eta|$. 
 \label{defn:cAinfty2}
\end{defn}
\begin{rem}
The sign of the cyclicity relation above is presented for instance in
\cite{Poli:higher}. 
\end{rem}
\begin{defn}[Minimal $A_\infty$-category]
An $A_\infty$-category $\cC$ 
is called {\em minimal} if $m_1=0$. 
 \label{defn:minimal}
\end{defn}

The suspension $s(\cC)$ of an $A_\infty$-category $\cC$ 
is defined by the shift 
\begin{equation}\label{suspension-cat}
 s: \cC(a,b)\to s(\cC(a,b))=:s(\cC)(a,b)
\end{equation}
for any $a,b\in\Ob(\cC)=\Ob(s(\cC))$.
Here, as a $\Z$-graded vector space,
we mean $(s(\cC))^r(a,b)=\cC^{r+1}(a,b)$. 
We denote the induced $A_\infty$-structure in $s(\cC)$ by
$\{\bar{m}_n\}_{n\ge 1}$. 
Explicitly, $\bar{m}_n$ is defined by 
\begin{equation}\label{m-suspension}
 \begin{split}
 & \bar{m}_n (s(w_{12}),\dots,s(w_{n(n+1)})) 
 =(-1)^\star s\cdot \bar{m}_n(w_{12},\dots,w_{n(n+1)}), \\
 & \star = {(n-1)|s(w_{12})|+(n-2)|s(w_{23})|+ \cdots +2|s(w_{(n-2)(n-1)})|
+ |s(w_{(n-1)n})|} , 
 \end{split}
\end{equation}
where the degree $|\bar{m}_n|$ becomes one for all $n\ge 1$. 
Let us denote $s(w_{i(i+1)})=:\bar{w}_{i(i+1)}\in s(\cC)(a_i,a_{i+1})$. 
Then, the $A_\infty$-relations
(\ref{Ainfty2}) turn out to be
\begin{equation}
 \begin{split}
&0=\sum_{k+l=n+1}\sum_{j=0}^{k-1}
(-1)^\star \\
&\ \bar{m}_k(\bar{w}_{12},\cdots,\bar{w}_{j(j+1)},\bar{m}_l(\bar{w}_{(j+1)(j+2)},\cdots,\bar{w}_{(j+l)(j+l+1)}),
 \bar{w}_{(j+l+1)(j+l+2)},\cdots,\bar{w}_{n(n+1)})\ ,
 \end{split}
 \label{Ainfty2}
\end{equation}
where $\star=|\bar{w}_{12}|+\cdots+|\bar{w}_{j(j+1)}|$. 
\begin{defn}[$A_\infty$-functor]
Given two $A_\infty$-categories $\cC$, $\cC'$, 
an {\em $A_\infty$-functor} 
$\bar{\f}:=\{f;\bar{f}_1,\bar{f}_2,\dots\}: s(\cC)\to s(\cC')$ is 
a map $f:\Ob(s(\cC))\to\Ob(s(\cC'))$ of objects with 
degree preserving multilinear maps 
\begin{equation*}
 \bar{f}_k: s(\cC)(a_1,a_2)\otimes\cdots\otimes s(\cC)(a_k,a_{k+1})
 \to s(\cC')(f(a_1),f(a_{k+1}))
\end{equation*}
satisfying 
\begin{equation*}
  \sum_{i\ge 1}\sum_{j_1+\cdots + j_i=n} \bar{m}'_i
  (\bar{f}_{j_1}\otimes\cdots\otimes \bar{f}_{j_i})
  = \sum_{k+l=n+1}\sum_{j=0}^{k-1}
  \bar{f}_k
 (\underbrace{id\otimes\cdots\otimes id}_j\otimes\bar{m}_l
 \otimes id\otimes\cdots\otimes id)  
\end{equation*}
for $n=1,2,\dots$. 
An $A_\infty$-functor $\bar{\f}$ with $\bar{f}_2=0,\ \bar{f}_3=0,\dots$
is called {\em linear}. 
On the other hand, 
if $f:\Ob(s(\cC))\to\Ob(s(\cC'))$ is a bijection 
and $\bar{f}_1: s(\cC)(a,b)\to s(\cC')(f(a),f(b))$ 
induces an isomorphism between the 
cohomologies for any $a,b\in\Ob(s(\cC))$, 
we call the $A_\infty$-functor an {\em $A_\infty$-quasi-isomorphism} functor. 
In particular, it is an {\em $A_\infty$-isomorphism} functor
if $\bar{f}_1$ is itself an isomorphism. 
\end{defn}

For a cyclic $A_\infty$-category $(\cC,\eta)$,
after the suspension $s:\cC\to s(\cC)$, 
the inner product $s(\eta)=:\omega$ in $s(\cC)$ is given by 
$\omega=\eta(s^{-1}\ ,s^{-1}\ )$, or more explicitly, 
\begin{equation*}
 \omega(\bar{w}_{ab},\bar{w}_{ba})=(-1)^{|\bar{w}_{ab}|}\eta(w_{ab},w_{ba}), 
\end{equation*}
where the cyclicity condition (\ref{cyclic-cd}) 
turns out to be \cite{hk:nchms}
\begin{equation*}
 \begin{split}
 & \omega(\bar{m}_n(\bar{w}_{12},\dots,\bar{w}_{n(n+1)}),\bar{w}_{(n+1)1})
 =(-1)^\star
 \omega(\bar{m}_n(\bar{w}_{23},\dots,\bar{w}_{(n+1)1}),\bar{w}_{12}),  \\
 & \star= (|\bar{w}_{23}|+\cdots +|\bar{w}_{(n+1)1}|)|\bar{w}_{12}|, 
 \end{split}
\end{equation*}
for homogeneous elements $\bar{w}_{i(i+1)}\in s(\cC)(a_i,a_{i+1})$, 
$i=1,\dots,n+1$ (with the identification $i+(n+1)=i$). 
\begin{defn}[Cyclic $A_\infty$-functor]
For two cyclic $A_\infty$-categories $\cC$ and $\cC'$ 
with the inner products $\eta$ and $\eta'$, respectively, 
we call an $A_\infty$-functor $\bar{f}:s(\cC)\to s(\cC')$ {\em cyclic} when 
\begin{equation}
  \omega'(\bar{f}_1(\bar{w}_{ab}),\bar{f}_1(\bar{w}_{ba}))=
  \omega(\bar{w}_{ab},\bar{w}_{ba}) 
 \label{omegacF1}
\end{equation}
and, for fixed $n\ge 3$, 
\begin{equation}
 \sum_{k,l\ge 1,\ k+l=n}
\omega'(\bar{f}_k(\bar{w}_{12},\dots,\bar{w}_{k(k+1)}),
 \bar{f}_l(\bar{w}_{(k+1)(k+2)},\dots,\bar{w}_{n(n+1)}))=0 
 \label{omegacF2}
\end{equation}
holds, where $\omega=s(\eta)$ and $\omega'=s(\eta')$. 
 \label{defn:cAinftyfunc}
\end{defn}

 \subsection{Triangulated $A_\infty$-category $\Tw(\cC)$}
\label{ssec:triAinfty}

In this subsection, 
we briefly
recall the construction of triangulated category $\Tri(\cC)$ from 
a given strictly unital $A_\infty$-category $\cC$
proposed in \cite{kont:hms}. 
This triangulated category $\Tri(\cC)$
is obtained as the zero-th cohomology category of the $A_\infty$-category
$\Tw(\cC)$ of one-sided twisted complexes in $\cC$. 
We first extend $\cC$ to the additive $A_\infty$-category $\ti{\cC}$
generated by $\cC$.
An object ${\bf a}\in\ti{\cC}$ is 
a formal finite direct sum of formally shifted objects
\begin{equation*}
 {\bf a}=a_{1}[r_1]\oplus \cdots \oplus a_{l}[r_l] . 
\end{equation*}
The $A_\infty$-structure $\{m_n\}$ is extended to the one $\{\ti{m}_n\}$ in $\ti{\cC}$. 
There exists a few natural choices of $\{\ti{m}_n\}$ which differ in signs.
In this paper, we set 
\begin{equation}\label{tim}
  \ti{m}_n(w'_{12},\dots,w'_{n(n+1)}):=
  (-1)^{nr_1} (m_n(w_{12},\dots,w_{n(n+1)}))' 
\end{equation}
for $w_{ij}\in\cC(a_i,a_j)$, and $w'$ implies the morphism corresponding to $w$ 
after the shifts $a_i \to a_i[r_i]$ by $r_i\in\Z$. 
This $\ti{m}_n$ is $\ti{m}_n^{(2)}$ in \cite{hk:Ainfty-sign}. 
This is because
we compare $\cC$ with a DG category in discussing the homological mirror
symmetry. 
The result $\ti{\cC}$ is again a strictly unital $A_\infty$-category. 
The corresponding $A_\infty$-products in $s(\ti{\cC})$ 
are also denoted by the same symbol $\ti{m}_n$. 

A twisted complex $({\bf a}, {\bf w})$ in $\ti{\cC}$ 
is an object 
\begin{equation*}
 {\bf a}=a_{1}[r_1]\oplus \cdots \oplus a_{l}[r_l] 
\end{equation*}
in $\ti{\cC}$ 
with an element ${\bf w}\in\ti s(\ti{\cC})^0({\bf a},{\bf a})$
satisfying the $A_\infty$-Maurer-Cartan equation
\begin{equation*}
  \ti{m}_1({\bf w})+\ti{m}_2({\bf w},{\bf w})
  + \ti{m}_3({\bf w},{\bf w},{\bf w}) +\cdots =0 . 
\end{equation*}
A twisted complex $({\bf a}, {\bf w})$ is called {\em one-sided}
if the $l\times l$ matrix
\begin{equation*}
  {\bf w} \in
  \bp
  s\ti{\cC}^0(a_{1}[r_1],a_{1}[r_1]) & \cdots & s\ti{\cC}^0(a_{1}[r_1],a_{l}[r_l])\\
  \vdots & \ddots & \vdots \\ 
  s\ti{\cC}^0(a_{l}[r_l],a_{1}[r_1]) & \cdots & s\ti{\cC}^0(a_{l}[r_l],a_{l}[r_l])
  \ep
\end{equation*}
is upper triangular and in addition that the diagonal entries are trivial. 
We denote by $\Tw(\cC)$ the category consisting of
one-sided twisted complexes in $\cC$, where 
the spaces of morphisms are defined by
$\Tw(\cC)(({\bf a},{\bf w}_a),({\bf b},{\bf w}_b)):=
\ti{\cC}({\bf a},{\bf b})$. An $A_\infty$-structure $\m^{\Tw}$ is
then given by 
\begin{equation*}
 \begin{split}
 &  m_n^{\Tw}({\bf w}_{12},\dots,{\bf w}_{n(n+1)}) \\
     & =\sum_{k_1,\dots,k_{n+1}\ge 0}
     \ti{m}_{n+k_1+\cdots +k_{n+1}}
     (({\bf w}_1)^{k_1},{\bf w}_{12},({\bf w}_2)^{k_2},\dots, {\bf w}_{n(n+1)},({\bf w}_{n+1})^{k_{n+1}})  
 \end{split}
\end{equation*}
for $({\bf a}_i,{\bf w}_i)\in\Tw(\cC)$ and
${\bf w}_{i(i+1)}\in s\Tw(\cC)(({\bf a}_i,{\bf w}_i),({\bf a}_{i+1},{\bf w}_{i+1}))$, $i=1,\dots, n+1$. 
Thus, $\Tw(\cC)$ again forms an $A_\infty$-category.

The zero-th cohomology of this $A_\infty$-category $\Tw(\cC)$
then forms a triangulated category \cite{kont:hms}, which we denote by
$H^0(\Tw(\cC)):=\Tri(\cC)$.  
The shift functor $T:\Tri(\cC)\to\Tri(\cC)$
is induced naturally from the additive automorphism
$T:\ti{\cC}\to\ti{\cC}$ such that $T(a)=a[1]$. See \cite{hk:Ainfty-sign}. 
Exact triangles in $\Tri(\cC)$ are defined as sequences which are
isomorphic to 
\begin{equation*}
  \cdots T^{-1}({\bf c},{\bf w}_c)\overset{{\bf w}}\to ({\bf a},{\bf w}_a)\to
  C({\bf w})\to ({\bf c}, {\bf w}_c)\to\cdots 
\end{equation*}
given by the {\em mapping cones} 
\begin{equation*}
  C({\bf w}):=\left({\bf c}\oplus {\bf a},
  \left( \bps {\bf w}_c & \bar{\bf w} \\ 0 & {\bf w}_a \eps\right)\right)
\end{equation*}
of $m_1^{\Tw(\cC)}$-closed 
morphisms ${\bf w}\in \Tw(\cC)(T^{-1}({\bf c},{\bf w}_c), ({\bf a}, {\bf w}_a))$
of degree zero, where 
${\bf w}'$ is the induced degree zero morphism in 
$s\Tw(\cC)(({\bf c},{\bf w}_c), ({\bf a},{\bf w}_a))$. 
For more details, see \cite{hk:Ainfty-sign} and references therein.

When $\cC$ is a strictly unital cyclic $A_\infty$-category, 
by direct calculations we can check the followings: 
\begin{prop}
For a strictly unital cyclic $A_\infty$-category $(\cC,\eta)$
of degree $|\eta|$, 
we set 
\begin{equation}\label{tieta}
 \ti\eta(w_{ab}',w'_{ba}):=(-1)^{(|\eta|+1)r_a}\eta(w_{ab},w_{ba}) 
\end{equation}
for any $w_{ab}$ and $w_{ba}$, where $w'_{**}$ is the morphism
corresponding to $w_{**}$ after the shifts $a\to a[r_a]$, $b\to b[r_b]$.

Then, this extends to a cyclic structure in $\ti\cC$.
\qed\end{prop}
Then, by the construction of $\Tw(\cC)$, it is clear that: 
\begin{cor}
$(\Tw(\cC),\ti\eta)$ forms a strictly unital cyclic $A_\infty$-category. 
\qed\end{cor}
Hence, this $\ti\eta$ induces the perfect pairings of the Serre duality
in the zero-th cohomology $H^0(\Tw(\cC))=\Tri(\cC)$.

 \section{Fukaya $A_\infty$ category of a two-torus}
\label{sec:Fukaya}
 
In this section, we construct 
the Fukaya category $\cC$ of a two-torus $T^2$ explicitly
as a strictly unital minimal $A_\infty$-category. 
As the objects, we treat Lagrangians consisting of {\em geodesic} cycles
in $T^2$ since
any nongeodesic cycle turns out be isomorphic to a geodesic one
in $\cC$. We define these objects in subsection \ref{ssec:objects}. 
The space of morphisms are also easy to define,
which is done in subsection \ref{ssec:morphisms}. 
What is complicated is the $A_\infty$-structure.
We spend a few subsections to preparation, and 
then define the $A_\infty$-structure
in subsection \ref{ssec:Fukaya}.

\subsection{The objects}
\label{ssec:objects}

Let us consider a two-torus $T^2$ 
whose covering space is $\ti{T}^2=\R^2$ 
with coordinates $(x,y)\in\R^2$. We have $\pi_{xy}\wt{T^2}=T^2$, 
$\pi_{xy}:=\pi_x\pi_y=\pi_y\pi_x$, 
where $\pi_x$ and $\pi_y$ are the projections 
associated with the identifications 
$x\sim x+1$ and $y\sim y+1$, respectively. 

We fix a flat complexified symplectic structure on $T^2$ and denote it by
\begin{equation*}
 \rho\cdot dx\wedge dy
\end{equation*}
with a complex number $\rho$ whose imaginary part is positive. 
Note that the standard symplectic structure corresponds to the case
$\rho= \ii$.

Let $p$ and $q$ be relatively prime integers 
such that $q\ge 0$, where $p=1$ if $q=0$. 
Consider a geodesic cycle $\pi_{xy}(L)\in T^2$, 
\begin{equation*}
 L : q y=p x+\alpha\ ,\qquad \alpha\in\R\ ,
\end{equation*}
with a number $\beta\in\R$. 
Then, $\pi^{-1}_{xy}\pi_{xy}(L)\subset\R^2$ is a copy of 
lines $q y=p x+\alpha+c$, $c\in\Z$. 
Each of the lines is described as 
\begin{equation*}
 (x,y)=(q l,p l)+ (x_0,y_0), \qquad l\in\R, 
\end{equation*}
with a fixed point $(x_0,y_0)\in\pi^{-1}_{xy}\pi_{xy}(L)$, 
and $\beta$ is regarded as a flat connection 
\begin{equation*}
 \left(\fpartial{l}-2\pi \ii\beta\right) dl
\end{equation*}
of a line bundle over the geodesic cycle. 
Here, note that the parameter $l$ is chosen so that 
the length of the geodesic cycle is one. 
We denote by $\Lag$ the set of all objects as above. 
Namely, an object $a\in\Ob(\cC)=\Lag$ is a quadruple
$a=(q_a,p_a,\alpha_a,\beta_a)$, 
where $(q_a,p_a)$ are relatively prime integers 
and $\alpha_a,\beta_a\in\R$. 
Thus, an object $a$ is associated with a geodesic
cycle $\pi_{xy}(L_a)$ 
and a flat connection on a line bundle over $\pi_{xy}(L_a)$
determined by $\beta_a$. 
Since $T^2$ is two-dimensional, any such (geodesic) cycle forms a
Lagrangian. 

For each object $a=(p_a,q_a,\alpha_a,\beta_a)\in\Lag$, 
we assign a number 
\begin{equation}
 -\frac{\pi}{2}<\phi_a:=\nat\left(\frac{p_a}{q_a}\right)\le\frac{\pi}{2} . 
\end{equation}

For given two objects $a,b\in\Lag$, one has 
$\pi^{-1}_{xy}\pi_{xy}(L_a)=\pi^{-1}_{xy}\pi_{xy}(L_b)$ 
if and only if $\phi_a=\phi_b$ and $\alpha_b-\alpha_a\in\Z$. 
We say $a$ is {\em isomorphic} to $b$ and denote it by $a\simeq b$ 
if and only if 
\begin{equation*}
 \phi_a=\phi_b,\qquad 
 \alpha_b-\alpha_a\in\Z,\qquad \beta_b-\beta_a\in\Z\ . 
\end{equation*}

 \subsection{The morphisms}
\label{ssec:morphisms}
 
Hereafter,
for any $a,b\in\Lag$, 
we often denote the space of morphisms by 
\begin{equation*}
  \cC(a,b)=:V_{ab} , 
\end{equation*}
which is set to be the following $\Z$-graded vector space over $\C$ 
of degree zero and one only. 

If $\phi_a\neq\phi_b$,
the corresponding geodesic cycles are transversal to each other,
and then $V_{ab}$ is by definition generated by 
the bases which are identified with
the intersection points of $(L_a,L_b)$ in $\T^2$. 
We denote $p_{ab}:=q_ap_b-p_aq_b$. 
We see that there exists $|p_{ab}|$ intersection
points of $\pi_{xy}(L_a)$ and $\pi_{xy}(L_b)$. 
We label these intersection points by $v_{ab}^j$, 
$j\in\Z/p_{ab}\Z$, and denote the corresponding base of $V_{ab}$ 
by the same symbol $v_{ab}^j$. 
Their degrees are defined as 
\begin{equation*}
 \begin{cases}
   |v_{ab}^j|=0 &  \phi_a<\phi_b, \\
   |v_{ab}^j|=1 & \phi_a>\phi_b
 \end{cases}
\end{equation*}
for any $j$. Thus, one has 
\begin{equation}
 \begin{cases}
  V_{ab}^0\simeq (\C)^{|p_{ab}|},\qquad V_{ab}^1=0 & \phi_a<\phi_b \\
  V_{ab}^0=0,\qquad V_{ab}^1\simeq (\C)^{|p_{ab}|} & \phi_a>\phi_b . 
 \end{cases}
\end{equation}

In case $\phi_{a}=\phi_b$, 
the corresponding geodesic cycles are not transversal to each other. 
One has 
$\pi_{xy}^{-1}\pi_{xy}(L_a)\cap\pi_{xy}^{-1}\pi_{xy}(L_b)=\emptyset$ 
if $a\not\simeq b$, and 
$\pi_{xy}^{-1}\pi_{xy}(L_a)=\pi_{xy}^{-1}\pi_{xy}(L_b)$
if $a\simeq b$. 
According to this, 
we set 
\begin{equation*}
 V_{ab}^0=0,\qquad V_{ab}^1=0
\end{equation*}
if $a\not\simeq b$, and 
\begin{equation*}
 V_{ab}^0\simeq\C,\qquad V_{ab}^1\simeq\C
\end{equation*}
if $a\simeq b$. 
We denote by $v_{ab}$ the base of $V_{ab}^0$ or 
$V_{ab}^1$. 

Let us define a degree minus one nondegenerate symmetric 
inner product $\eta:V_{ab}\otimes V_{ba}\to\C$, 
$a,b\in\Ob(\cC)$. 
For $a,b\in\Ob(\cC)$ such that $\phi_a\ne\phi_b$, 
we set 
\begin{equation*}
 \eta(v_{ab},v_{ba})=
 \begin{cases}
  1 & v_{ab}=v_{ba} \ \ \text{in $T^2$}\\
  0 &  \text{otherwise}. 
 \end{cases}
\end{equation*}
For $a\simeq b\in\Ob(\cC)$, 
we set 
\begin{equation*}
 \eta(v_{ab},v_{ba})=1 
\end{equation*}
if $|v_{ab}|+|v_{ba}|=1$ and zero otherwise. 
Then, 
for any base $v_{ab}\in V_{ab}$, 
the pairing $\eta$ defines a dual base which 
we denote by $(v_{ab})^*\in V_{ba}^{1-|v_{ab}|}$. 

This $\eta$ turns out to define a cyclic structure of degree one
in the $A_\infty$-category $\cC$ which we will define
in subsection \ref{ssec:Fukaya}.

 \subsection{Clockwise convex (CC-) polygons}
\label{ssec:CCpolygon}
 
\def \aa{{\frak a}}

In order to define the $A_\infty$-structure in $\cC$, 
we introduce the notions of CC-collections, CC-sequences, 
CC-polygons and so on. 
Let $\aa:=(a_1,\dots,a_{n+1})$, $a_1,\dots,a_{n+1}\in\Lag$, 
be a collection of objects such that 
$\phi_{a_i}<\phi_{a_{i+1}}$ for two of $i=1,\dots, n+1$ 
and $\phi_{a_i}>\phi_{a_{i+1}}$ for other $i$, 
where we identify $a_{j+(n+1)}$ with $a_j$. 
We call such a collection $\aa$ a {\em CC-collection}. 
We also call the isomorphism classes 
$[\aa]:=([a_1],\dots, [a_{n+1}]$) also a CC-collection. 
Note that the corresponding base $v_{a_ia_{i+1}}$
is of degree zero (resp. one) when 
$\phi_{a_i}<\phi_{a_{i+1}}$ (resp. $\phi_{a_i}>\phi_{a_{i+1}}$). 

Given a CC collection $\aa$, consider a sequence 
$\vec{v}:=(v_{a_1a_2},\dots,v_{a_na_{n+1}},v_{a_{n+1}a_1})$ 
of intersection points in $T^2$, 
where $v_{a_ia_{i+1}}\in\pi_{xy}(L_{a_i}\cap L_{a_{i+1}})$. 
We call such $\vec{v}$ a {\em CC-sequence}. 
We set the degrees by $|v_{a_ia_{i+1}}|=0$ if $\phi_{a_i}<\phi_{a_{i+1}}$ 
and $|v_{a_ia_{i+1}}|=1$ if $\phi_{a_i}>\phi_{a_{i+1}}$. 
These are exactly the degrees of the bases of morphisms
corresponding to the points $v_{**}$. 

For a given CC-sequence $\vec{v}$, 
let $CC'(\vec{v})$ be a subset of 
\begin{equation*}
 \{\ti{v}:=(\ti{v}_{a_1a_2},\dots,\ti{v}_{a_na_{n+1}}, v_{a_{n+1}a_1})
 \in (\pi_{xy}^{-1}({v}_{a_1a_2}),\dots,\pi_{xy}^{-1}(\ti{v}_{a_na_{n+1}}), 
\pi_{xy}^{-1}(v_{a_{n+1}a_1}))\}
\end{equation*}
consisting of all elements $\ti{v}$ satisfying the followings: 
\begin{itemize}
 \item  for each $i=1,\dots,n+1$, if 
$\ti{v}_{a_{i-1}a_i}\ne\ti{v}_{a_ia_{i+1}}$, then the geodesic interval 
between $\ti{v}_{a_{i-1}a_i}$ and $\ti{v}_{a_ia_{i+1}}$ is included in 
$\pi_{xy}^{-1}\pi_{xy}(L_{a_i})$, 
 \item $\ti{v}_{a_{n+1}a_1}=v_{a_{n+1}a_1}\in\R^2$, where 
we fixed an inclusion of the fundamental domain of $T^2$ 
to the covering space $\R^2$ and denoted the image of 
${v}_{a_{n+1}a_1}$ also by ${v}_{a_{n+1}a_1}$ itself. 
\end{itemize}
We call any element $\ti{v}\in CC'(\vec{v})$ 
a {\em CC-semi-polygon} in the covering space $\R^2$. 

Let us express $\ti{v}\in CC'(\vec{v})$ as the form
\begin{equation*}
 \ti{v}=(v_1,\dots,v_1,v_2,\dots,v_2,\dots\dots, v_n,\dots,v_n) ,  
\end{equation*}
where $\{v_1,\dots,v_n\}$, $n\in\Z_{> 0}$, are points in $\R^2$ 
such that $v_i\ne v_{i+1}$ for $i=1,\dots, n-1$. 
In this expression, 
we call $\ti{v}$ a {\em CC-point} if $n=1$, 
a {\em CC-line} if $n=2$ or $n=3$ and $v_1=v_3$, 
and a {\em CC-polygon} 
otherwise. 
\footnote{This CC-polygon in $\R^2$ is 
that defined in \cite{hk:Ainftyplane}. }
We denote by $CC(\vec{v})\subset CC'(\vec{v})$ the subset
consisting of a CC-point and CC-polygons. 
For generic $\vec{v}$, there does not exist any CC-line 
and in this case $CC(\vec{v})= CC'(\vec{v})$.

For a CC-line or a CC-polygon 
$\ti{v}=(v_1,\dots,v_1,v_2,\dots,v_2,\dots\dots, v_n,\dots,v_n)$, 
let $v_i$ and $v_j$, $i<j$, be two degree zero vertices. 
Then, the pair $(x(v_i), x(v_j))$ gives left/right or right/left 
extrema of $x(v_k)$, $k\in\{1,\dots, n+1\}$.

\begin{figure}[h]
 \begin{minipage}[c]{75mm}{
\begin{center}
\unitlength 0.1in
\begin{picture}( 26.7500, 20.1500)( 19.7500,-29.1500)
%
\special{pn 8}%
\special{pa 4200 2400}%
\special{pa 3900 2600}%
\special{pa 3550 2700}%
\special{pa 3100 2750}%
\special{pa 2800 2700}%
\special{pa 2500 2550}%
\special{pa 2300 2300}%
\special{pa 2200 2000}%
\special{pa 2250 1600}%
\special{pa 2400 1300}%
\special{pa 2600 1100}%
\special{pa 2900 1000}%
\special{pa 3300 950}%
\special{pa 3700 1000}%
\special{pa 4100 1100}%
\special{pa 4250 1200}%
\special{pa 4400 1400}%
\special{pa 4500 1800}%
\special{pa 4450 2100}%
\special{pa 4350 2300}%
\special{pa 4200 2400}%
\special{pa 4200 2400}%
\special{pa 4200 2400}%
\special{fp}%
%
\special{pn 20}%
\special{sh 1}%
\special{ar 4200 2400 10 10 0  6.28318530717959E+0000}%
\special{sh 1}%
\special{ar 4200 2400 10 10 0  6.28318530717959E+0000}%
%
\special{pn 20}%
\special{sh 1}%
\special{ar 3900 2600 10 10 0  6.28318530717959E+0000}%
\special{sh 1}%
\special{ar 3900 2600 10 10 0  6.28318530717959E+0000}%
%
\special{pn 20}%
\special{sh 1}%
\special{ar 3500 2700 10 10 0  6.28318530717959E+0000}%
\special{sh 1}%
\special{ar 3500 2700 10 10 0  6.28318530717959E+0000}%
%
\special{pn 20}%
\special{sh 1}%
\special{ar 3100 2750 10 10 0  6.28318530717959E+0000}%
\special{sh 1}%
\special{ar 3100 2750 10 10 0  6.28318530717959E+0000}%
%
\special{pn 20}%
\special{sh 1}%
\special{ar 2800 2700 10 10 0  6.28318530717959E+0000}%
\special{sh 1}%
\special{ar 2800 2700 10 10 0  6.28318530717959E+0000}%
%
\special{pn 20}%
\special{sh 1}%
\special{ar 2500 2550 10 10 0  6.28318530717959E+0000}%
\special{sh 1}%
\special{ar 2500 2550 10 10 0  6.28318530717959E+0000}%
%
\special{pn 20}%
\special{sh 1}%
\special{ar 2300 2300 10 10 0  6.28318530717959E+0000}%
\special{sh 1}%
\special{ar 2300 2300 10 10 0  6.28318530717959E+0000}%
%
\special{pn 20}%
\special{sh 1}%
\special{ar 2200 2000 10 10 0  6.28318530717959E+0000}%
\special{sh 1}%
\special{ar 2200 2000 10 10 0  6.28318530717959E+0000}%
%
\special{pn 20}%
\special{sh 1}%
\special{ar 2250 1600 10 10 0  6.28318530717959E+0000}%
\special{sh 1}%
\special{ar 2250 1600 10 10 0  6.28318530717959E+0000}%
%
\special{pn 20}%
\special{sh 1}%
\special{ar 2400 1300 10 10 0  6.28318530717959E+0000}%
\special{sh 1}%
\special{ar 2400 1300 10 10 0  6.28318530717959E+0000}%
%
\special{pn 20}%
\special{sh 1}%
\special{ar 2600 1100 10 10 0  6.28318530717959E+0000}%
\special{sh 1}%
\special{ar 2600 1100 10 10 0  6.28318530717959E+0000}%
%
\special{pn 20}%
\special{sh 1}%
\special{ar 2900 1000 10 10 0  6.28318530717959E+0000}%
\special{sh 1}%
\special{ar 2900 1000 10 10 0  6.28318530717959E+0000}%
%
\special{pn 20}%
\special{sh 1}%
\special{ar 3300 950 10 10 0  6.28318530717959E+0000}%
\special{sh 1}%
\special{ar 3300 950 10 10 0  6.28318530717959E+0000}%
%
\special{pn 20}%
\special{sh 1}%
\special{ar 3700 1000 10 10 0  6.28318530717959E+0000}%
\special{sh 1}%
\special{ar 3700 1000 10 10 0  6.28318530717959E+0000}%
%
\special{pn 20}%
\special{sh 1}%
\special{ar 4100 1100 10 10 0  6.28318530717959E+0000}%
\special{sh 1}%
\special{ar 4100 1100 10 10 0  6.28318530717959E+0000}%
%
\special{pn 20}%
\special{sh 1}%
\special{ar 4250 1200 10 10 0  6.28318530717959E+0000}%
\special{sh 1}%
\special{ar 4250 1200 10 10 0  6.28318530717959E+0000}%
%
\special{pn 20}%
\special{sh 1}%
\special{ar 4400 1400 10 10 0  6.28318530717959E+0000}%
\special{sh 1}%
\special{ar 4400 1400 10 10 0  6.28318530717959E+0000}%
%
\special{pn 20}%
\special{sh 1}%
\special{ar 4500 1800 10 10 0  6.28318530717959E+0000}%
\special{sh 1}%
\special{ar 4500 1800 10 10 0  6.28318530717959E+0000}%
%
\special{pn 20}%
\special{sh 1}%
\special{ar 4450 2100 10 10 0  6.28318530717959E+0000}%
\special{sh 1}%
\special{ar 4450 2100 10 10 0  6.28318530717959E+0000}%
%
\special{pn 20}%
\special{sh 1}%
\special{ar 4350 2300 10 10 0  6.28318530717959E+0000}%
\special{sh 1}%
\special{ar 4350 2300 10 10 0  6.28318530717959E+0000}%
\put(22.5000,-20.5000){\makebox(0,0)[lb]{$v_i$}}%
\put(44.5000,-18.5000){\makebox(0,0)[rb]{$v_j$}}%
\put(35.0000,-26.0000){\makebox(0,0){$v_1$}}%
\put(31.0000,-26.5000){\makebox(0,0){$v_2$}}%
\put(28.0000,-26.0000){\makebox(0,0){$v_3$}}%
%
\special{pn 8}%
\special{pa 3500 2450}%
\special{pa 3468 2456}%
\special{pa 3438 2462}%
\special{pa 3406 2468}%
\special{pa 3374 2474}%
\special{pa 3342 2478}%
\special{pa 3310 2484}%
\special{pa 3278 2488}%
\special{pa 3248 2492}%
\special{pa 3216 2494}%
\special{pa 3184 2498}%
\special{pa 3152 2500}%
\special{pa 3120 2500}%
\special{pa 3088 2500}%
\special{pa 3056 2500}%
\special{pa 3024 2498}%
\special{pa 2992 2496}%
\special{pa 2960 2492}%
\special{pa 2928 2486}%
\special{pa 2896 2480}%
\special{pa 2864 2472}%
\special{pa 2834 2462}%
\special{pa 2804 2452}%
\special{pa 2774 2440}%
\special{pa 2744 2426}%
\special{pa 2716 2410}%
\special{pa 2688 2392}%
\special{pa 2660 2374}%
\special{pa 2634 2354}%
\special{pa 2610 2332}%
\special{pa 2586 2310}%
\special{pa 2564 2286}%
\special{pa 2542 2262}%
\special{pa 2522 2236}%
\special{pa 2506 2208}%
\special{pa 2490 2180}%
\special{pa 2476 2152}%
\special{pa 2464 2122}%
\special{pa 2454 2092}%
\special{pa 2446 2060}%
\special{pa 2438 2028}%
\special{pa 2432 1996}%
\special{pa 2430 1964}%
\special{pa 2426 1932}%
\special{pa 2426 1898}%
\special{pa 2426 1866}%
\special{pa 2430 1832}%
\special{pa 2432 1798}%
\special{pa 2438 1766}%
\special{pa 2444 1734}%
\special{pa 2450 1700}%
\special{pa 2458 1668}%
\special{pa 2468 1638}%
\special{pa 2480 1606}%
\special{pa 2492 1576}%
\special{pa 2504 1546}%
\special{pa 2518 1518}%
\special{pa 2534 1490}%
\special{pa 2552 1464}%
\special{pa 2570 1438}%
\special{pa 2590 1412}%
\special{pa 2612 1388}%
\special{pa 2634 1366}%
\special{pa 2658 1344}%
\special{pa 2682 1324}%
\special{pa 2708 1304}%
\special{pa 2736 1288}%
\special{pa 2764 1270}%
\special{pa 2792 1254}%
\special{pa 2820 1240}%
\special{pa 2850 1228}%
\special{pa 2880 1216}%
\special{pa 2910 1206}%
\special{pa 2942 1196}%
\special{pa 2972 1186}%
\special{pa 3004 1180}%
\special{pa 3036 1172}%
\special{pa 3068 1168}%
\special{pa 3100 1162}%
\special{pa 3132 1158}%
\special{pa 3164 1156}%
\special{pa 3196 1154}%
\special{pa 3228 1152}%
\special{pa 3260 1150}%
\special{pa 3292 1150}%
\special{pa 3326 1150}%
\special{pa 3358 1152}%
\special{pa 3390 1154}%
\special{pa 3420 1156}%
\special{pa 3452 1158}%
\special{pa 3484 1162}%
\special{pa 3516 1166}%
\special{pa 3548 1170}%
\special{pa 3580 1176}%
\special{pa 3610 1182}%
\special{pa 3642 1188}%
\special{pa 3674 1194}%
\special{pa 3706 1202}%
\special{pa 3738 1210}%
\special{pa 3768 1218}%
\special{pa 3800 1228}%
\special{pa 3830 1238}%
\special{pa 3860 1250}%
\special{pa 3890 1264}%
\special{pa 3918 1278}%
\special{pa 3944 1296}%
\special{pa 3970 1316}%
\special{pa 3994 1336}%
\special{pa 4016 1360}%
\special{pa 4038 1384}%
\special{pa 4056 1410}%
\special{pa 4074 1436}%
\special{pa 4092 1462}%
\special{pa 4106 1490}%
\special{pa 4120 1520}%
\special{pa 4132 1550}%
\special{pa 4144 1580}%
\special{pa 4154 1612}%
\special{pa 4162 1644}%
\special{pa 4170 1676}%
\special{pa 4176 1708}%
\special{pa 4180 1740}%
\special{pa 4184 1774}%
\special{pa 4184 1806}%
\special{pa 4184 1840}%
\special{pa 4182 1872}%
\special{pa 4178 1904}%
\special{pa 4172 1936}%
\special{pa 4164 1966}%
\special{pa 4152 1996}%
\special{pa 4140 2024}%
\special{pa 4124 2054}%
\special{pa 4108 2080}%
\special{pa 4088 2106}%
\special{pa 4068 2132}%
\special{pa 4046 2156}%
\special{pa 4024 2178}%
\special{pa 4000 2202}%
\special{pa 3976 2222}%
\special{pa 3950 2242}%
\special{pa 3926 2262}%
\special{pa 3900 2282}%
\special{pa 3874 2300}%
\special{pa 3848 2318}%
\special{pa 3822 2336}%
\special{pa 3800 2350}%
\special{sp}%
%
\special{pn 8}%
\special{pa 3822 2336}%
\special{pa 3800 2350}%
\special{fp}%
\special{sh 1}%
\special{pa 3800 2350}%
\special{pa 3868 2330}%
\special{pa 3844 2320}%
\special{pa 3844 2296}%
\special{pa 3800 2350}%
\special{fp}%
%
\special{pn 8}%
\special{pa 2050 2850}%
\special{pa 4650 2850}%
\special{fp}%
\special{sh 1}%
\special{pa 4650 2850}%
\special{pa 4584 2830}%
\special{pa 4598 2850}%
\special{pa 4584 2870}%
\special{pa 4650 2850}%
\special{fp}%
%
\special{pn 8}%
\special{pa 2200 900}%
\special{pa 2200 2850}%
\special{dt 0.045}%
%
\special{pn 8}%
\special{pa 4500 900}%
\special{pa 4500 2850}%
\special{dt 0.045}%
\put(22.0000,-30.0000){\makebox(0,0){$x(v_i)$}}%
\put(45.0000,-30.0000){\makebox(0,0){$x(v_j)$}}%
\end{picture}%

\vspace*{0.4cm}
 (a)
\end{center}}
 \end{minipage}
\quad 
 \begin{minipage}[c]{75mm}{
\begin{center}
\unitlength 0.1in
\begin{picture}( 26.7500, 20.1500)( 19.7500,-29.1500)
%
\special{pn 8}%
\special{pa 4200 2400}%
\special{pa 3900 2600}%
\special{pa 3550 2700}%
\special{pa 3100 2750}%
\special{pa 2800 2700}%
\special{pa 2500 2550}%
\special{pa 2300 2300}%
\special{pa 2200 2000}%
\special{pa 2250 1600}%
\special{pa 2400 1300}%
\special{pa 2600 1100}%
\special{pa 2900 1000}%
\special{pa 3300 950}%
\special{pa 3700 1000}%
\special{pa 4100 1100}%
\special{pa 4250 1200}%
\special{pa 4400 1400}%
\special{pa 4500 1800}%
\special{pa 4450 2100}%
\special{pa 4350 2300}%
\special{pa 4200 2400}%
\special{pa 4200 2400}%
\special{pa 4200 2400}%
\special{fp}%
%
\special{pn 20}%
\special{sh 1}%
\special{ar 4200 2400 10 10 0  6.28318530717959E+0000}%
\special{sh 1}%
\special{ar 4200 2400 10 10 0  6.28318530717959E+0000}%
%
\special{pn 20}%
\special{sh 1}%
\special{ar 3900 2600 10 10 0  6.28318530717959E+0000}%
\special{sh 1}%
\special{ar 3900 2600 10 10 0  6.28318530717959E+0000}%
%
\special{pn 20}%
\special{sh 1}%
\special{ar 3500 2700 10 10 0  6.28318530717959E+0000}%
\special{sh 1}%
\special{ar 3500 2700 10 10 0  6.28318530717959E+0000}%
%
\special{pn 20}%
\special{sh 1}%
\special{ar 3100 2750 10 10 0  6.28318530717959E+0000}%
\special{sh 1}%
\special{ar 3100 2750 10 10 0  6.28318530717959E+0000}%
%
\special{pn 20}%
\special{sh 1}%
\special{ar 2800 2700 10 10 0  6.28318530717959E+0000}%
\special{sh 1}%
\special{ar 2800 2700 10 10 0  6.28318530717959E+0000}%
%
\special{pn 20}%
\special{sh 1}%
\special{ar 2500 2550 10 10 0  6.28318530717959E+0000}%
\special{sh 1}%
\special{ar 2500 2550 10 10 0  6.28318530717959E+0000}%
%
\special{pn 20}%
\special{sh 1}%
\special{ar 2300 2300 10 10 0  6.28318530717959E+0000}%
\special{sh 1}%
\special{ar 2300 2300 10 10 0  6.28318530717959E+0000}%
%
\special{pn 20}%
\special{sh 1}%
\special{ar 2200 2000 10 10 0  6.28318530717959E+0000}%
\special{sh 1}%
\special{ar 2200 2000 10 10 0  6.28318530717959E+0000}%
%
\special{pn 20}%
\special{sh 1}%
\special{ar 2250 1600 10 10 0  6.28318530717959E+0000}%
\special{sh 1}%
\special{ar 2250 1600 10 10 0  6.28318530717959E+0000}%
%
\special{pn 20}%
\special{sh 1}%
\special{ar 2400 1300 10 10 0  6.28318530717959E+0000}%
\special{sh 1}%
\special{ar 2400 1300 10 10 0  6.28318530717959E+0000}%
%
\special{pn 20}%
\special{sh 1}%
\special{ar 2600 1100 10 10 0  6.28318530717959E+0000}%
\special{sh 1}%
\special{ar 2600 1100 10 10 0  6.28318530717959E+0000}%
%
\special{pn 20}%
\special{sh 1}%
\special{ar 2900 1000 10 10 0  6.28318530717959E+0000}%
\special{sh 1}%
\special{ar 2900 1000 10 10 0  6.28318530717959E+0000}%
%
\special{pn 20}%
\special{sh 1}%
\special{ar 3300 950 10 10 0  6.28318530717959E+0000}%
\special{sh 1}%
\special{ar 3300 950 10 10 0  6.28318530717959E+0000}%
%
\special{pn 20}%
\special{sh 1}%
\special{ar 3700 1000 10 10 0  6.28318530717959E+0000}%
\special{sh 1}%
\special{ar 3700 1000 10 10 0  6.28318530717959E+0000}%
%
\special{pn 20}%
\special{sh 1}%
\special{ar 4100 1100 10 10 0  6.28318530717959E+0000}%
\special{sh 1}%
\special{ar 4100 1100 10 10 0  6.28318530717959E+0000}%
%
\special{pn 20}%
\special{sh 1}%
\special{ar 4250 1200 10 10 0  6.28318530717959E+0000}%
\special{sh 1}%
\special{ar 4250 1200 10 10 0  6.28318530717959E+0000}%
%
\special{pn 20}%
\special{sh 1}%
\special{ar 4400 1400 10 10 0  6.28318530717959E+0000}%
\special{sh 1}%
\special{ar 4400 1400 10 10 0  6.28318530717959E+0000}%
%
\special{pn 20}%
\special{sh 1}%
\special{ar 4500 1800 10 10 0  6.28318530717959E+0000}%
\special{sh 1}%
\special{ar 4500 1800 10 10 0  6.28318530717959E+0000}%
%
\special{pn 20}%
\special{sh 1}%
\special{ar 4450 2100 10 10 0  6.28318530717959E+0000}%
\special{sh 1}%
\special{ar 4450 2100 10 10 0  6.28318530717959E+0000}%
%
\special{pn 20}%
\special{sh 1}%
\special{ar 4350 2300 10 10 0  6.28318530717959E+0000}%
\special{sh 1}%
\special{ar 4350 2300 10 10 0  6.28318530717959E+0000}%
\put(22.5000,-20.5000){\makebox(0,0)[lb]{$v_j$}}%
\put(44.5000,-18.5000){\makebox(0,0)[rb]{$v_i$}}%
\put(37.0000,-11.0000){\makebox(0,0){$v_1$}}%
\put(40.2000,-11.8000){\makebox(0,0){$v_2$}}%
\put(41.7000,-12.8000){\makebox(0,0){$v_3$}}%
%
\special{pn 8}%
\special{pa 3600 1260}%
\special{pa 3634 1266}%
\special{pa 3666 1272}%
\special{pa 3698 1278}%
\special{pa 3730 1284}%
\special{pa 3762 1292}%
\special{pa 3792 1302}%
\special{pa 3822 1312}%
\special{pa 3852 1324}%
\special{pa 3880 1338}%
\special{pa 3906 1354}%
\special{pa 3932 1372}%
\special{pa 3956 1392}%
\special{pa 3980 1414}%
\special{pa 4002 1438}%
\special{pa 4022 1464}%
\special{pa 4040 1490}%
\special{pa 4058 1518}%
\special{pa 4074 1548}%
\special{pa 4090 1576}%
\special{pa 4104 1606}%
\special{pa 4116 1636}%
\special{pa 4126 1668}%
\special{pa 4134 1698}%
\special{pa 4142 1730}%
\special{pa 4146 1762}%
\special{pa 4150 1792}%
\special{pa 4152 1826}%
\special{pa 4150 1858}%
\special{pa 4148 1890}%
\special{pa 4142 1922}%
\special{pa 4136 1954}%
\special{pa 4126 1984}%
\special{pa 4116 2014}%
\special{pa 4104 2044}%
\special{pa 4088 2074}%
\special{pa 4072 2102}%
\special{pa 4056 2128}%
\special{pa 4036 2154}%
\special{pa 4016 2180}%
\special{pa 3994 2204}%
\special{pa 3970 2226}%
\special{pa 3946 2248}%
\special{pa 3922 2270}%
\special{pa 3896 2290}%
\special{pa 3868 2308}%
\special{pa 3842 2326}%
\special{pa 3814 2344}%
\special{pa 3786 2358}%
\special{pa 3756 2374}%
\special{pa 3728 2388}%
\special{pa 3698 2400}%
\special{pa 3668 2412}%
\special{pa 3638 2424}%
\special{pa 3608 2434}%
\special{pa 3576 2444}%
\special{pa 3546 2452}%
\special{pa 3514 2460}%
\special{pa 3482 2466}%
\special{pa 3452 2472}%
\special{pa 3420 2478}%
\special{pa 3388 2484}%
\special{pa 3356 2488}%
\special{pa 3324 2492}%
\special{pa 3290 2494}%
\special{pa 3258 2498}%
\special{pa 3226 2500}%
\special{pa 3194 2500}%
\special{pa 3162 2502}%
\special{pa 3128 2502}%
\special{pa 3096 2502}%
\special{pa 3064 2502}%
\special{pa 3032 2500}%
\special{pa 3000 2496}%
\special{pa 2968 2492}%
\special{pa 2936 2488}%
\special{pa 2904 2482}%
\special{pa 2874 2474}%
\special{pa 2842 2466}%
\special{pa 2812 2454}%
\special{pa 2782 2442}%
\special{pa 2752 2430}%
\special{pa 2724 2414}%
\special{pa 2696 2398}%
\special{pa 2668 2380}%
\special{pa 2642 2360}%
\special{pa 2616 2340}%
\special{pa 2592 2316}%
\special{pa 2570 2294}%
\special{pa 2548 2270}%
\special{pa 2528 2244}%
\special{pa 2510 2216}%
\special{pa 2494 2188}%
\special{pa 2480 2160}%
\special{pa 2468 2130}%
\special{pa 2456 2100}%
\special{pa 2448 2070}%
\special{pa 2440 2038}%
\special{pa 2434 2006}%
\special{pa 2430 1974}%
\special{pa 2428 1940}%
\special{pa 2426 1908}%
\special{pa 2426 1874}%
\special{pa 2428 1840}%
\special{pa 2432 1808}%
\special{pa 2436 1774}%
\special{pa 2442 1742}%
\special{pa 2448 1710}%
\special{pa 2456 1678}%
\special{pa 2466 1646}%
\special{pa 2476 1616}%
\special{pa 2488 1586}%
\special{pa 2500 1556}%
\special{pa 2514 1526}%
\special{pa 2530 1498}%
\special{pa 2548 1472}%
\special{pa 2566 1446}%
\special{pa 2584 1420}%
\special{pa 2606 1396}%
\special{pa 2628 1372}%
\special{pa 2650 1350}%
\special{pa 2676 1328}%
\special{pa 2700 1308}%
\special{pa 2728 1290}%
\special{pa 2756 1274}%
\special{pa 2784 1258}%
\special{pa 2812 1246}%
\special{pa 2842 1234}%
\special{pa 2874 1224}%
\special{pa 2904 1216}%
\special{pa 2936 1210}%
\special{pa 2968 1204}%
\special{pa 3000 1200}%
\special{pa 3032 1198}%
\special{pa 3064 1196}%
\special{pa 3098 1196}%
\special{pa 3130 1196}%
\special{pa 3162 1198}%
\special{pa 3194 1200}%
\special{pa 3226 1202}%
\special{pa 3258 1206}%
\special{pa 3290 1210}%
\special{pa 3300 1210}%
\special{sp}%
%
\special{pn 8}%
\special{pa 3290 1210}%
\special{pa 3300 1210}%
\special{fp}%
\special{sh 1}%
\special{pa 3300 1210}%
\special{pa 3236 1184}%
\special{pa 3248 1206}%
\special{pa 3232 1224}%
\special{pa 3300 1210}%
\special{fp}%
%
\special{pn 8}%
\special{pa 2200 900}%
\special{pa 2200 2850}%
\special{dt 0.045}%
%
\special{pn 8}%
\special{pa 4500 900}%
\special{pa 4500 2850}%
\special{dt 0.045}%
\put(22.0000,-30.0000){\makebox(0,0){$x(v_j)$}}%
\put(45.0000,-30.0000){\makebox(0,0){$x(v_i)$}}%
%
\special{pn 8}%
\special{pa 2050 2850}%
\special{pa 4650 2850}%
\special{fp}%
\special{sh 1}%
\special{pa 4650 2850}%
\special{pa 4584 2830}%
\special{pa 4598 2850}%
\special{pa 4584 2870}%
\special{pa 4650 2850}%
\special{fp}%
\end{picture}%

\vspace*{0.4cm}
 (b)
\end{center}}
 \end{minipage}
 \caption{CC-polygons $\vec{v}$ 
with (a): $\sigma(\vec{v})=-1$ and (b): $\sigma(\vec{v})=+1$. }
\label{fig:polygon-sign}
\end{figure}
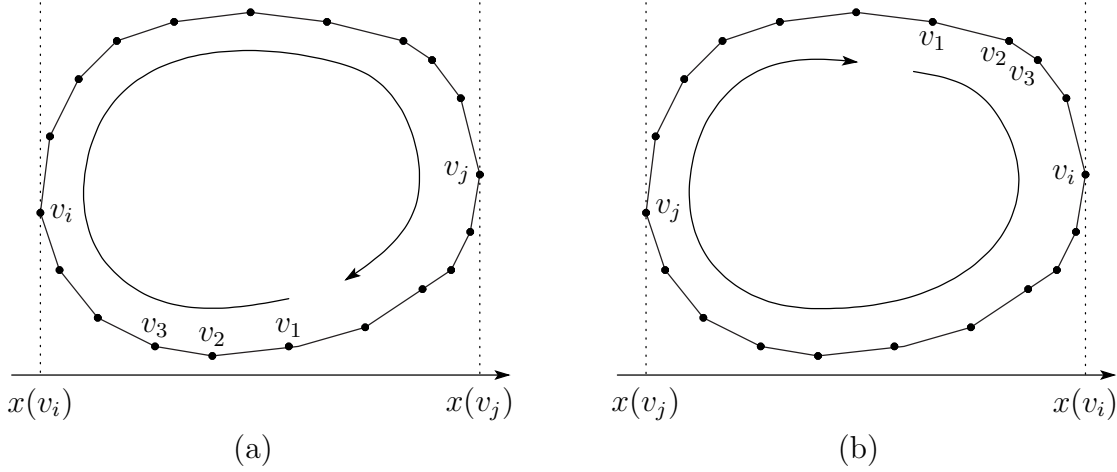

The {\em sign $\sigma(\ti{v})$} of a CC semi-polygon $\ti{v}$ 
is then defined by 
\begin{equation}\label{sigma}
 \sigma(\ti{v}):=
 \begin{cases}
   -1 & x(v_i) < x(v_j) \\
   +1 & x(v_j) < x(v_i) 
 \end{cases}  
\end{equation}
if $\ti{v}$ is a CC-line or a CC-polygon 
(see Figure \ref{fig:polygon-sign}).  

For a CC-point $\ti{v}$, we set $\sigma(\ti{v})=+1$. 

We define the {\em multiplicity} $d_i$ of $v_i$, $i=1,\dots, n$, 
as follows. 
\begin{itemize}
 \item 
If $\ti{v}$ forms a CC-line or a CC-polygon, then, 
the multiplicity of $(v_i)^{\otimes 1+d_i}$ is $d_i$ for 
$v_i$ of degree zero and 
the multiplicity of $(v_i)^{\otimes d_i}$ is $d_i$ for 
$v_i$ of degree one. 

Note that, when $\ti{v}$ forms a line, then, $n=2$ and 
both $v_1$ and $v_2$ are degree one. 

We attach to $v_i$ a number 
\begin{equation*}
\begin{cases}
 D_i:=\ov{(d_i)!} & |v_i|=1 \\
 D_i:=\ov{2^{d_i}(d_i)!} & |v_i|=0. 
\end{cases}
\end{equation*}

\item If $\ti{v}$ forms a point, then $n=1$ and 
we set the multiplicity of $(v_1)^{\otimes 2+d_i}$ as $d_1$ 
and then 
\begin{equation*}
 D_1:=
 \begin{cases}
 1 &  d_1=1 \\
 0 & \text{otherwise}.
\end{cases}
\end{equation*}
\end{itemize}
For any element $\ti{v}\in CC'(\vec{v})$ in this description, 
we set $V(\ti{v}):=D_1\cdots D_n$ and 
denote by $Area(\ti{v})$ the area of the CC-semi-polygon. 
The symplectic area is then $-\ii\rho\cdot Area(\ti{v})$. 
Note that $Area(\ti{v})=0$ if $\ti{v}$ forms a CC-point or a CC-line.

 \subsection{Constants $F(\vec{v},\beta)$ and $F(\vec{v},{\frak B})$}
\label{ssec:generalizedF}
 
In this subsection,
we define constants $F(\vec{v},\beta)$ and their generalizations
$F(\vec{v};{\frak B})$. 
These will appear 
as the structure constants, i.e.,
the coefficients of the $A_\infty$-products 
in the Fukaya $A_\infty$-category $\cC$ of a two-torus. 
\begin{defn}\label{defn:F}
Given a CC-sequence $\vec{v}=(v_{a_1a_2},\dots,v_{a_na_{n+1}},v_{a_{n+1}a_1})$ 
and \\
$\beta:=(\beta_{a_1},\dots,\beta_{a_{n+1}})\in\R^{n+1}$, 
define 
\begin{equation}
 \begin{split}
& F(\vec{v};\beta)
:= \sum_{\ti{v}\in CC(\vec{v})}
 (\sigma(\ti{v}))^{n}\ 
V(\ti{v})\exp(2\pi\ii\rho\cdot Area(\ti{v}))
\exp{(2\pi\ii\sum_{i=1}^{n+1}\beta_{a_{i}}(\ti{v}))}, \\
 & \beta_{a_i}(\ti{v}):=\beta_{a_i}\cdot\Delta_i l(\ti{v}),
 \end{split}
 \label{fukaya}
\end{equation}
where, $\Delta_i l(\ti{v})\in\R$ is determined by 
\begin{equation*}
  (x,y)(\ti{v}_{a_ia_{i+1}})-(x,y)(\ti{v}_{a_{i-1}a_i})
 =\Delta_i l(\ti{v})(p_{a_i},q_{a_i}). 
\end{equation*}
This implies 
$\Delta_i l(\ti{v})=(x(\ti{v}_{a_ia_{i+1}})-x(\ti{v}_{a_{i-1}a_i}))/q_{a_i}$
if $q_{a_i}\ne 0$ and 
$\Delta_i l(\ti{v})=(y(\ti{v}_{a_ia_{i+1}})-y(\ti{v}_{a_{i-1}a_i}))/p_{a_i}$
if $p_{a_i}\ne 0$. 
\end{defn}
Note that 
$\Delta_i l(\ti{v})=0$ 
if $\ti{v}_{a_ia_{i+1}}=\ti{v}_{a_{i-1}a_i}$.

Next, we define a generalization of $F(\vec{v},\beta)$. 

Given a CC-collection $[\aa]=([a_1],\dots,[a_{n+1}])$, 
let us fix $(b_1,\dots,b_{n+1})\in (\Z_{\ge 0})^{n+1}$ and 
consider a collection $a_i^{\lambda_i}$, 
$\lambda_i=0,\dots, b_i$, of objects which are isomorphic to $a_i$ 
for each $i=1,\dots, n+1$. 
Denote $a_i^{\lambda_i}=(p_i, q_i, \alpha_i^{\lambda_i},\beta_i^{\lambda_i})$. 
By definition, $\beta_i^{\lambda_i'}-\beta_i^{\lambda_i}\in\Z$ 
for $\lambda_i,\lambda_i'\in\{0,\dots,b_i\}$. 
We denote 
\begin{equation*}
{\frak B}:=(\beta_{a_1}^0,\dots,\beta_{a_1}^{b_1};\cdots; 
\beta_{a_{n+1}}^{b_{n+1}},\dots,\beta_{a_{n+1}}^{b_{n+1}}) . 
\end{equation*}
Given an element $\ti{v}=(\ti{v}_{a_1a_2},\dots,\ti{v}_{a_na_{n+1}},v_{a_{n+1}a_1})
\in CC(\vec{v})$ with $\vec{v}$ a CC-sequence of $\aa$,
define 
$\exp{(2\pi\ii\int{\frak B}_i(\ti{v}))}$, where $i=1,\dots,n+1$, by 
\begin{equation*}
 \int_{(l(\ti{v}_{a_{i-1}a_i})=l_i^0,\dots, 
l_i^{\lambda_i},\dots,
l_i^{b_i}=l(\ti{v}_{a_ia_{i+1}}))}
 \exp{\left(2\pi\ii (\beta_i^0(l_i^1-l_i^0)+\cdots 
 +\beta_i^{b_i}(l_i^{b_i}-l_i^{b_i-1}) )\right)} . 
\end{equation*}
Here the integral $\int_{(l(\ti{v}_{a_{i-1}a_i})=l_i^0,\dots,
l_i^{\lambda_i},\dots,l_i^{b_i}=l(\ti{v}_{a_ia_{i+1}}))}$ indicates 
\begin{equation*}
 \begin{cases}
  \int_{l_i^0\le l_i^1\le\cdots\le l_i^{b_i-1}
\le l_i^{b_i}}
 dl_i^1\cdots dl_i^{b_i-1} & 
l(\ti{v}_{a_{i-1}a_i})\le l(\ti{v}_{a_ia_{i+1}}), \\
  \int_{l_i^{b_i}\le l_i^{b_i-1}
 \le\cdots\le l_i^{1}
\le l_i^{0}}
 dl_i^{b_i-1}\cdots dl_i^{1} & 
 l(\ti{v}_{a_{i-1}a_i})\ge l(\ti{v}_{a_ia_{i+1}}) 
 \end{cases}
\end{equation*}
for $l_i^0=l(\ti{v}_{a_{i-1}a_i})\in\R$ and 
$l_i^{b_i}=l(\ti{v}_{a_ia_{i+1}})\in\R$ 
satisfying 
\begin{equation*}
 (x,y)(\ti{v}_{a_ia_{i+1}})-(x,y)(\ti{v}_{a_{i-1}a_i})
 =(l_i^{b_i}-l_i^0)(p_{a_i},q_{a_i}). 
\end{equation*}
We set 
\begin{equation}
 F(\vec{v};{\frak B}):=\sum_{\ti{v}\in CC(\vec{v})}
 (\sigma(\ti{v}))^{n+b}\ 
V(\ti{v})
\exp(2\pi\ii\rho\cdot Area(\ti{v}))
 \prod_i\exp{(2\pi\ii\int{\frak B}_i(\ti{v}))}  
\end{equation}
where $b:=b_1+\cdots +b_n$. 

By direct calculations, one obtains the followings. 
\begin{lem}\label{lem:multi-beta}
{\em (i)}\ \ When $b=0$, 
one has $F(\vec{v};{\frak B})=F(\vec{v};\beta)$ by
the identification  
$\beta_i^0=\beta_{a_i}$ for any $i=1,\dots, n+1$. 

{\em (ii)}\ \ For generic ${\frak B}$, one has 
\begin{equation*}
\exp{(2\pi\ii\int{\frak B}_i(\ti{v}))}
 = \left( \ov{2\pi\ii}\right)^{b_i}
\sum_{j\in\{0,\dots,b_i\}}
\left(
\frac{\exp{\left(2\pi\ii\beta_i^j\cdot\Delta_i l(\ti{v})\right)}}{\prod_{k\in\{0,\dots,b_i\}\backslash \{j\}}(\beta_i^j-\beta_i^k)}
\right). 
\end{equation*}

{\em (iii)}\ \
If $\beta_i^0=\cdots=\beta_i^{b_i}=:\beta_i$, then 
\begin{equation*}
 \begin{split}
\exp{(2\pi\ii\int{\frak B}_i(\ti{v}))}
 & = \frac{|\Delta_i l(\ti{v})|^{b_i}}{(b_i)!}
   \exp{\left(2\pi\ii\beta_{a_i}\cdot \Delta_i l(\ti{v})\right)} \\
 & =\left(\ov{(b_i)!}\left(\frac{\sigma(i)\sigma(\ti{v})}{2\pi\ii}\fpartial{\beta_{a_i}}\right)^{b_i}
 \right)
   \exp{\left(2\pi\ii\beta_{a_i}\cdot\Delta_i l(\ti{v})\right)} . 
\end{split}
\end{equation*}
where $\sigma(i)=-1$ if $k \le i \le l$ 
with two elements $v_{a_{k-1}a_k}$ and $v_{a_{l}a_{l+1}}$ 
of degree zero 
and $\sigma(i)=+1$ otherwise. 

Thus, if $\beta_i^0=\cdots=\beta_i^{b_i}=:\beta_{a_i}$ for any $i$, then 
\begin{equation*}\label{Fformula}
 F(\vec{v};{\frak B})=  \prod_{i=1}^{n+1}
 \left(\ov{(b_i)!}\left(\frac{\sigma(i)}{2\pi\ii}\fpartial{\beta_{a_i}}
\right)^{b_i}\right) F(\vec{v};\beta) . 
\end{equation*}
\qed\end{lem}

 \subsection{The structure constants $\cF$}
\label{ssec:str-const}

We will see in subsection \ref{ssec:Fukaya} that 
the generalizations $F(\vec{v};{\frak B})$ 
will be the structure constants of 
generic $A_\infty$-products in $\cC$. 
Recall that they are associated with CC-polygons. 
However, an $A_\infty$-product may sometimes has contributions
from CC-lines. 
In this subsection,
we introduce such a correction $R(\vec{v};{\frak B})$ 
so that the sum $F(\vec{v};{\frak B})+ R(\vec{v};{\frak B})$
will be the structure constants of any
$A_\infty$-products in $\cC$.

{}For a given CC-sequence $\vec{v}$, 
let $L(\vec{v};a_i,a_j)$, $2\le i<j\le n+1$, be the set of all CC-lines 
$\ti{v}\in CC'(\vec{v})\backslash CC(\vec{v})$ such that 
\begin{itemize}
 \item $\beta_j-\beta_i\in\Z$, 
 \item 
$\Delta_k l(\ti{v})\ne 0$
for $k=i$ and $k=j$. 
\end{itemize}
The set $L(\vec{v};a_i,a_j)$ can be nonempty 
only if $a_i\simeq a_j$. 

Let us prepare the following functions
\begin{equation}\label{r^d}
 r^d(x)=  \sum_{n\in\Z\backslash \{0\}} \ov{(2\pi\ii n)^d} 
 \exp(-2\pi\ii n x)  
\end{equation}
in variable $x$, where $d=1,2,\dots$ .
These are periodic, that is,
$r^d(x)=r^d(x+1)$ holds.
By the inverse Fourier transformation, one sees that 
$r^d$ is a polynomial of degree $d$ in $(0,1)$
such that $\int_0^1 r^d(x) dx =0$. 
Though $r^1$ is discontinuous at $x\in\Z$, 
the value is zero, $r^1(\Z)=0$, since it is a periodic odd function. 
For $d>1$, the function $r^d$ becomes continuous 
since $r^{d+1} = (r^d)'$ holds. 
\footnote{The treatment of the value $r^d(\Z)$ seems incorrect
in \cite{hk:hpt-hms}. }

Then, we set 
\begin{equation*}
 R(\vec{v};{\frak B}):=
 -\sum_{\{(i<j)| L(\vec{v};a_i,a_j)\ne\emptyset\}} \delta_{b_i+b_j,b}
\frac{b!}{b_i!b_j!} 4V(\ti{v}) r^{b+1}(\Delta_i l(\ti{v}))   . 
\end{equation*}
In the right hand side, $\ti{v}$ are those belonging to $L(\vec{v};a_i,a_j)$. 
We see that the result does not depend on the choice of such an $\ti{v}$. 

The sum 
\begin{equation}\label{str-const}
 \cF(\vec{v},{\frak B}):=
 F(\vec{v};{\frak B}) + R(\vec{v};{\frak B})\  
\end{equation}
will be the structure constants of the $A_\infty$-products in $\cC$.

 \subsection{The Fukaya $A_\infty$-category of a two-torus}
\label{ssec:Fukaya}
 
\begin{defn}[Fukaya $A_\infty$-category $\cC$ of a a flat symplectic two-torus]
The set of objects is $\Ob(\cC):=\Lag$ (Subsection \ref{ssec:objects}). 
For any $a,b\in\Lag$, the space 
$\cC(a,b)=:V_{ab}$ of morphisms is a graded vector space
over $\C$ 
of degree zero and one as already defined in subsection \ref{ssec:morphisms}. 
In particular, if $\phi_a\neq\phi_b$, the space $V_{ab}$ is generated by 
the basis $\{v_{ab}^j\}_{j\in\Z/p_{ab}\Z}$
associated to the intersection points of $\pi_{xy}(L_a,L_b)$ in $T^2$.

The degree minus one nondegenerate symmetric 
inner product $\eta:V_{ab}\otimes V_{ba}\to\C$ is also already defined.
Recall that,
for a morphism $v_{ab}\in V_{ab}$,
the dual base $(v_{ab})^*\in V_{ba}$
satisfies $\eta(v_{ab},(v_{ab})^*)=1$.

Now, we define a minimal cyclic $A_\infty$-structure in $\cC$. 
For any $n\ge 2$ and $a_1,\dots, a_{n+1}\in\Ob(\cC)$, 
we define a collection $\{\varphi_{n+1}\}_{n\ge 2}$ of
multilinear maps 
\begin{equation*}
 \varphi_{n+1}: V_{a_1a_2}\otimes\cdots\otimes V_{a_{n}a_{n+1}}
 \otimes V_{a_{n+1}a_1}\to\C , 
\end{equation*}
of degree $1-n$ which are cyclic, i.e., which satisfy 
\begin{equation*}
 \begin{split}
 & \varphi_{n+1}(w_{a_1a_2},\dots,w_{a_na_{n+1}},w_{a_{n+1}a_1})
 =(-1)^\star
\varphi_{n+1}(w_{a_{n+1}a_1},w_{a_1a_2},\dots,w_{a_na_{n+1}}) \\
 & \quad \star = n+|w_{a_{n+1}a_1}|(|w_{a_1a_2}|+\cdots +|w_{a_na_{n+1}}|)
 \end{split}
 \end{equation*}
for homogeneous elements 
$w_{a_ia_{i+1}}\in V_{a_ia_{i+1}}$, $i=1,\dots,n$. 
Here, that the degree is $1-n$ implies that 
\begin{equation*}
  |w_{a_1a_2}|+\cdots +|w_{a_na_{n+1}}|+|w_{a_{n+1}a_1}|=n-1 . 
\end{equation*}

First,
if $a_i\simeq a_{i+1}$ and $|w_{a_ia_{i+1}}|=0$ for some $i\in\{1,\dots,n\}$, 
then we set
\begin{equation*}
  \varphi_{n+1}(w_{a_1a_2},\dots,w_{a_na_{n+1}},w_{a_{n+1}a_1})=0
\end{equation*}
for $n\ge 3$. 
For $n=2$, for any $a,a',b\in\Ob(\cC)$ such that $a\simeq a'$, 
define a cyclic tri-linear map 
$\varphi_3:V_{ab}\otimes V_{ba'}\otimes V_{a'a}^0\to\C$ by 
\begin{equation}\label{imply-identity}
 \varphi_3(v_{ab},v_{ba'},v_{a'a})=1 
\end{equation}
if and only if $\phi_a\ne\phi_b$ and $v_{ba'}=v_{ab}$ 
or $a\simeq b$ and $|v_{ab}|+|v_{ba}|=1$, and zero otherwise.
This determines all the cyclic tri-linear maps 
$\varphi_3:V_{a_1a_2}\otimes V_{a_2a_3}\otimes V_{a_3a_1}\to\C$ 
such that 
$a_i\simeq a_{i+1}$ and $|w_{a_ia_{i+1}}|=0$ for some $i\in\{1,\dots,n\}$.

Next, for $n\ge 2$ and a given 
$(b_1,\dots,b_{n+1})\in(\Z_{\ge 0})^{\otimes (n+1)}$, 
consider a collection of isomorphic objects $a_i^{\lambda_i}$, 
$\lambda_i\in\{0,\dots,b_i\}$ for each $i$ 
such that $\phi_{[a_i]}\ne\phi_{[a_{i+1}]}$ for any $i=1,\dots,n+1$. 
Denote  $a_i^0=:a_i$, $a_i^{b_i}=:a_i'$, $V_{a_i\cdots a_i'}^1
:=V_{a_i^0a_i^1}^1\otimes V_{a_i^1a_i^2}^1\otimes\cdots\otimes 
V_{a_i^{b_i-1}a_i^{b_i}}^1$ and 
$v_{a_i\cdots a_i'}:=
v_{a_i^0a_i^1}\otimes\cdots\otimes v_{a_i^{b_i-1}a_i^{b_i}}
\in V_{a_i^0\cdots a_i^{b_i}}^1$. 

Assume $\phi_{a_i}\ne \phi_{a_{i+1}}$ for $i=1,\dots,n$, 
and we set a cyclic map 
\begin{equation*}
\varphi_{n+1+b}: 
 V^1_{a_1\cdots a_1'}\otimes V_{a_1'a_2^0}\otimes
V^1_{a_2\cdots a_2'}\otimes\cdots
\otimes V_{a_{n+1}\cdots a_{n+1}'}^1
\otimes V_{a_{n+1}'a_1}\to\C
\end{equation*}
by 
\begin{equation}\label{varphi-def}
 \varphi_{n+1+b}(v_{a_1\cdots a_1'},v_{a_1'a_2},
v_{a_2\cdots a_2'},\dots, 
v_{a_{n+1}\cdots a_{n+1}'}, 
v_{a_{n+1}'a_1}) :=\cF(\vec{v},{\frak B}) 
\end{equation}
for $\vec{v}=(v_{a_1'a_2},\dots,v_{a_n'a_{n+1}},v_{a_{n+1}'a_1})$ 
if $\{[a_1],\dots,[a_{n+1}]\}$ is a CC-collection, and zero otherwise, 
where $\cF(\vec{v},{\frak B})$ is the structure constant given in
(\ref{str-const}). 
One sees that the definition (\ref{varphi-def})
is at least compatible with the cyclicity. 
We extend this $\varphi_{n+1+b}$ so that it will be cyclic. 

The above data determine the collection
$\{\varphi_{n+1}\}_{n\ge 2}$ of 
all cyclic multilinear maps of degree $1-n$. 
Define multilinear maps $m_n$, $n\ge 1$, of degree $2-n$ by 
$m_1=0$ and 
\begin{equation*}
 m_n(w_{a_1a_2},\dots,w_{a_na_{n+1}}):=
 \sum_{v_{a_1a_{n+1}}\in\Z/p_{a_1a_{n+1}}\Z}
 \varphi(w_{a_1a_2},\dots,w_{a_na_{n+1}},(v_{a_1a_{n+1}})^*)
 \cdot v_{a_1a_{n+1}}\ 
\end{equation*}
for $n\ge 2$ and $w_{a_ia_{i+1}}\in V_{a_ia_{i+1}}$, $i=1,\dots,n$. 
 \label{defn:fukaya}
\end{defn}
\begin{thm}\label{thm:main}
This $\cC=(\Lag,V=\oplus_{a,b\in\Lag}V_{ab},\eta,\m)$ 
forms a strictly unital minimal cyclic $A_\infty$-category. 
\end{thm}
The proof is obtained by showing 
homological mirror symmetry for (non)commutative two-tori. 
We explain this in section \ref{sec:hms}. 
\begin{rem}
The equation (\ref{imply-identity}) 
indicates that 
$v_{aa'}\in V_{aa'}^0$ for $a\simeq a'$ gives isomorphisms 
\begin{equation*}
 V_{a'b}\overset{\sim}{\to} V_{ab},\qquad 
 V_{ba}\overset{\sim}{\to} V_{ba'},\qquad 
\end{equation*}
by $m_2(v_{aa'}, * )$ 
and $m_2(* , v_{aa'})$, 
respectively. 
In particular, for any $a\in\Lag$, $v_{aa}\in V_{aa}^0$ is the 
unit in this $A_\infty$-category. 
\end{rem}

 \subsection{An explicit basis of the space of morphisms}
\label{ssec:basis}
 
As mentioned in subsection \ref{ssec:morphisms}, 
in this subsection, we give 
a natural way of attaching $v_{ab}^j$, $j\in\Z/p_{ab}\Z$
for given two objects $a$ and $b$. 
The reader can skip this subsection, 
but we employ this description fully later
in subsection \ref{ssec:m2}.

For $a=(q_a,p_a,\alpha_a,\beta_a)$,
we first attach integers $r_a$ and $s_a$ so that
\begin{equation*}
  \bp
   q_a & s_a \\
   p_a & r_a
  \ep
\end{equation*} 
is an $SL(2;\Z)$ element. 
Similarly, we attach an $SL(2;\Z)$ matrix
$\left(\bps q_b & s_b\\ p_b & r_b\eps\right)$ for $b$. 
We further define 
\begin{equation*}
   \bp
   q_{ab} & s_{ab} \\
   p_{ab} & r_{ab}
  \ep
  :=
    \bp
   q_a & s_a \\
   p_a & r_a
  \ep^{-1}
    \bp
   q_b & s_b \\
   p_b & r_b
   \ep =
     \bp
   r_aq_b-s_ap_b & r_as_b-s_ar_b \\
   q_ap_b-p_aq_b & q_ar_b-p_as_b
  \ep ,
\end{equation*}
where recall that $|p_{ab}|=\dim V_{ab}$
if $p_{ab}\ne 0$ (subsection \ref{ssec:morphisms}).  

Then, we define $v_{ab}^j$ as the projection of the intersection point of 
\begin{equation*}
 L_a : q_a y=p_a x+\alpha_a\ ,\qquad 
 L'_b : q_b\left(y-p_a\frac{q_{ab}j}{p_{ab}}\right)
 =p_b\left(x-q_a\frac{q_{ab}j}{p_{ab}}\right)+\alpha_b
\end{equation*}
in $\R^2$ by $\pi_{xy}:\R^2\to T^2$. 
Though $q_{ab}=r_aq_b-s_ap_b$ includes $r_a$ and $s_a$, 
we see that different choices of $(r_a,s_a)$ correspond to
permutations of $\{j=0,1,\cdots,|p_{ab}|-1\}$.

For later convenience, we mention that the defining equation for $L'_b$
is rewritten for instance as 
\begin{equation*}
L'_b : q_b\left(y-p_a\frac{q_{ab}j-\alpha_b}{p_{ab}}\right)
 =p_b\left(x-q_a\frac{q_{ab}j-\alpha_b}{p_{ab}}\right)
\end{equation*}
or 
\begin{equation*}
 L'_b : q_by =p_b x -q_{ab}j +\alpha_b . 
\end{equation*}

We also note that, if we define $v_{ba}^{j'}$ in a similar way as above, 
then we have 
\begin{equation*}
 (x,y)(v_{ba}^{-q_{ab}j})=(x,y)(v_{ab}^j) , 
\end{equation*}
i.e., the dual base of $v_{ab}^j$ is $(v_{ab}^j)^*=v_{ba}^{-q_{ab}j}$.

 \subsection{Fukaya $A_\infty$-categories for noncommutative two-tori}
\label{ssec:ncFukaya}
 
As is discussed in \cite{hk:foliation}, 
the mirror dual of the category of holomorphic vector bundles 
over a noncommutative complex torus 
turns out to be the Fukaya category on the mirror dual symplectic 
torus {\em equipped with a Kronecker foliation structure}, 
where the noncommutative parameter $\theta$ corresponds to
the slope of the foliations. 
In this subsection, we briefly mention
what we should modified from $\cC$ in the previous subsection
in order to obtain the Fukaya category $\cC_\theta$ depending on $\theta$. 
See also \cite[Section 5]{hk:rims07nc},
where it is explained more explicitly. 

Geometrically, we may regard the commutative case $\cC=\cC_{\theta=0}$
as the one which has the foliation structure defined by $x=\text{const}$. 
For general $\theta$, we define $\cC_\theta$ so that
it is associated with the foliation structure given by
$x+\theta y=\text{const}$.
In particular, we treat the line defined by $x+\theta y =0$ 
as if it is the $y$-axis. 
We denote by $x^\theta:=x+\theta y$ the resulting $x$-coordinate. 
Then, we first replace $\phi_a$ in
subsection \ref{ssec:morphisms} by 
\begin{equation*}
  -\frac{\pi}{2} <
  \mu^\theta_a:=\tan^{-1}\left( \frac{p_a}{q_a+p_a\theta}\right)
  \le\frac{\pi}{2} , 
\end{equation*}
and define the degree of a morphism in $\cC_\theta$
in a similar way as in subsection \ref{ssec:morphisms}. 

If we take
\begin{equation*}
 0<\tan^{-1}\left(\ov{\theta}\right) < \pi,  
\end{equation*}
where $\tan^{-1}(1/ 0)=: \pi/2$, 
we can also replace $\phi_a$ by 
\begin{equation*}
  \tan^{-1}\left(\ov{\theta}\right) -\pi <
  \phi^\theta_a:=\tan^{-1}\left( \frac{p_a}{q_a}\right)
  \le \tan^{-1}\left(\ov{\theta}\right) . 
\end{equation*}
Looking this $\phi^\theta_a$ in $x^\theta y$-plane gives $\mu^\theta_a$. 
Thus, $\mu^\theta_a\le \mu^\theta_b$
(resp. $\mu^\theta_a\ge \mu^\theta_b$)
if and only if 
$\phi^\theta_a\le \phi^\theta_b$
(resp. $\phi^\theta_a\ge \phi^\theta_b$).

Next, in the definition of the sign $\sigma(\ti{v})$
in (\ref{sigma}), 
$x$ is replaced by $x^\theta$. 
These are the all changes we need, and we obtain $\cC_\theta$. 
For instance, the effect of the holonomy $\beta_{a_i}(\ti{v})$
in Definition \ref{defn:F} 
does not depend on $\theta$;
the effects of the modification of $\beta_{a_i}$
and that of $\Delta_il(\ti{v})$ by $\theta$ cancel with each other,
so we may not care about them. 

The result $\cC_\theta$ again forms a minimal cyclic
$A_\infty$-algebra. 
Actually, it does not depend on $\theta$ so much
in the sense that $\ti\cC_\theta$ is naturally linear $A_\infty$-isomorphic
to $\ti\cC_{\theta'}$ for any $\theta$ and $\theta'$. 
We explain this fact in subsection \ref{ssec:theta-indep}.

\section{The homological mirror symmetry}
\label{sec:hms}

The construction of the explicit $A_\infty$-structure in the
Fukaya category of a two-torus in Theorem \ref{thm:main}
is due to 
applying the homological perturbation theory (HPT) 
to a DG category of holomorphic vector bundles on the mirror dual
complex torus. 
The most technical part of this construction is already done in
\cite{hk:Ainftyplane}. 
In this subsection, we briefly explain this
with assuming various arguments in references cited below. 

For a fixed $\theta$, 
we start with the cyclic DG category $\cC_{DG,\theta}$ of
indecomposable finitely generated projective right
modules over a noncommutative two torus $\cA_\theta$ equipped with
constant curvature connections as in \cite[Section 5.2]{hk:rims07nc}. 
Each object $a\in\cC_{DG,\theta}$ is in one-to-one correspondence
with the data $(q_a,p_a,\alpha_a,\beta_a)$. 
In particular, 
$\cC_{DG,\theta=0}$ corresponds to the DG category
$DG_{\check{T}^2}$ of holomorphic vector bundles
discussed in \cite{kobayashi:2tori}. 
Actually, $q_a$ and $p_a$ are the rank and the first Chern class,
respectively, of the corresponding holomorphic vector bundle
in $DG_{\check{T}^2}$. 
The correspondence of morphisms between these two DG categories
can be obtained
by Fourier expanding the spaces of morphisms in $DG_{\check{T}^2}$. 
However, for $\theta=0$, we can not include $a$ such that $q_a=0$. 
So, we choose an irrational $\theta$, 
and show
\begin{equation*}
  \cC_\theta\simeq\cC_{DG,\theta} .
\end{equation*}
This use of $\theta$ is auxiliary 
since $\ti\cC_{DG,\theta}$ is independent of $\theta$. See also
subsection \ref{ssec:theta-indep}. 

We would like to apply the homological perturbation theory (HPT) 
to $\cC_{DG,\theta}$ to obtain the corresponding Fukaya category.  
However, since the Fukaya category $\cC$ which we should obtain
corresponds to a particular limit $\epsilon\to 0$ of $A_\infty$-categories
obtained by the HPT, 
we can not obtain $\cC_\theta$ in such a direct way. 
See Introduction of \cite{hk:Ainftyplane} or \cite[Subsection 4.E]{hk:fukayadeform}.
Thus, we need to replace $\cC_{DG,\theta}$ with
a DG quasi-isomorphic DG category $\cC'_{DG,\theta}$ 
which includes delta function one-forms and step functions
as morphisms.

This replacement is an analogue of the replacement of 
$\cC_{DR}({\frak F})$ with $\cC'_{DR}({\frak F})$ 
we did in \cite{hk:Ainftyplane},
but we further need to modify a little more as we explain below.

First, it is easy to derive the $F(\vec{v},B)$ part of the
structure constants. 
We call the corresponding $A_\infty$-products
{\em transversal products}. 
Namely, an $A_\infty$-product 
\begin{equation*}
  m_n:V_{a_1a_2}\otimes\cdots\otimes V_{a_{n}a_{n+1}}\to V_{a_1a_{n+1}}
\end{equation*}
in $\cC_\theta$ is transversal if and only if 
$a_i\not\simeq a_{i+1}$ ($i\in\Z/n\Z$) and
$L(\vec{v};a_i,a_j)=\emptyset$ for any $1\le i<j\le n+1$. 
Thus, $m_n$ is always a transversal product
if $a_i\not\simeq a_j$ for any $1\le i<j\le n+1$. 
In order to derive these transversal products,
we do not need any non-transversal morphism. 
Thus, we may set
$\cC'_{DG,\theta}$ just as an analog of $\cC'_{DR}({\frak F})$
and apply the HPT to $\cC'_{DG,\theta}$. 
Any transversal product is then obtained as a sum of
transversal products 
of the Fukaya category of $\R^2$ discussed in \cite{hk:Ainftyplane}, 
where the sum runs over CC-polygons $\ti{v}$ associated
to a CC-sequence $\vec{v}$.

In order to derive the full $A_\infty$-products,
we need to include the spaces of non-transversal morphisms and
construct their Hodge decompositions. 
First,  
the space $\cC'_{DG,\theta}(a,a)$ should include
$1_a^*$. We treat this 
as $1^*_a:=\int_0^1dx  (d\vartheta_x)$, where
$d\vartheta_x$ is the delta function one forms
concentrated on $x+ \Z$ 
which is similar to the $d\vartheta_x\in\cC'_{DR}({\frak F})(a,a)$
in \cite{hk:Ainftyplane} but is now periodic. 
Then, we can derive the $F(\vec{v},{\frak B})$ parts 
of the structure constants $\cF(\vec{v},{\frak B})$. 
Actually, they are obtained by integrating 
results for the Fukaya category of $\R^2$ in \cite{hk:Ainftyplane}
over $x\in [0,1]$ and summing over the corresponding CC-polygons. 

We further need the following extension to the strategy in \cite{hk:Ainftyplane}
which is slightly complicated.
Since the DG category $\cC'_{DG,\theta}$ should
be closed with respect to the composition, 
the space $\cC'_{DG,\theta}(a,b)$ should include
the periodic delta function one forms $d\vartheta_x$ with each $x$,
not only their integration $\int_0^1dx  (d\vartheta_x)$. 
Here, if $a\not\simeq b$, then
the cohomology of the space $\cC'_{DG,\theta}(a,b)$ is trivial. 
Thus, the choice of a homotopy operator,
which determines a Hodge decomposition, is unique. 
However, if $a\simeq b$, then both
$H^0(\cC'_{DG,\theta}(a,b))$ and
$H^1(\cC'_{DG,\theta}(a,b))$ should be both one-dimensional.
Correspondingly, 
the homotopy operator $h_{ab}$ should be slightly modified. 
Then, it turns out that 
$h_{ab}(d\vartheta_c)=r^1(x-c)$,
where $r^1$ is the function defined in
(\ref{r^d}). 
Higher terms are derived, for instance for $a=b$, 
as $h_{ab}(r^1(x)\cdot 1_a^*)=r^2(x)$, $h_{ab}(r^2(x)\cdot 1_a^*)=r^3(x)$, ... . 
Furthermore, even for a transversal pair $(a,b)$,
the space $\cC'_{DG,\theta}(a,b)$ need to include more morphisms
created by the products
$\cC'_{DG,\theta}(a,a)\times\cC'_{DG,\theta}(a,b)$ and
$\cC'_{DG,\theta}(a,b)\times\cC'_{DG,\theta}(b,b)$. 
We see that these $r^d$ create the terms $R^d(\vec{v};a_i,a_j)$ of the
structure constants. 

We construct the DG category $\cC'_{DG,\theta}$ and the Hodge decomposition
in this way. 
Though the obtained $A_\infty$-products in $\cC_\theta$ are expressed by
infinite sums, 
the HPT guarantees that they are well defined.
Furthermore, 
the DG category $\cC'_{DG,\theta}$ has a natural cyclic structure, 
and the Hodge decomposition above preserves the cyclicity
in the sense in \cite[Lemma 5.23]{hk:thesis}. 
It is easier to see that it preserves the strict units. 
Thus, we can conclude that what we constructed above are
a strictly unital cyclic $A_\infty$-category $\cC_\theta$ and 
a strictly unital cyclic $A_\infty$-quasi-isomorphism
\begin{equation*}
  \cC_\theta\to\cC_{DG,\theta} . 
\end{equation*}

 \section{Symmetries of the functions $F$}
\label{sec:auto}

Recall that, for a CC-sequence 
\begin{equation*}
  \vec{v}=(v_{a_1a_2},\dots,v_{a_na_{n+1}},v_{a_{n+1}a_1})
\end{equation*}
each element $v_{a_ia_{i+1}}\in T^2$ is an intersection point 
of $\pi_{xy}(L_{a_i})$ and $\pi_{xy}(L_{a_{i+1}})$ 
labeled by $j_{a_ia_{i+1}}\in\Z/p_{a_ia_{i+1}}\Z$. 
Thus, $\vec{v}$ is determined by 
$\phi:=\{\phi_{a_1},\dots,\phi_{a_{n+1}}\}$, 
$J:=\{j_{a_1a_2},\dots,j_{a_na_{n+1}},j_{a_{n+1}a_1}\}$, 
and $\alpha_{a_1},\dots,\alpha_{n+1}$. 
In this sense, 
$F(\vec{v};\beta)$ is regarded as a function, determined by $\phi$ and $J$, 
in variables 
\begin{equation*}
(\alpha,\beta):=
  (\alpha_{a_1},\dots,\alpha_{a_{n+1}};\beta_{a_1},\dots,\beta_{a_{n+1}}) .
\end{equation*}
For a general $\theta\ne 0$, the only change is to replace $\phi$ 
by $\phi^\theta:=\{\phi^\theta_{a_1},\dots,\phi^\theta_{a_{n+1}}\}$. 
In this sense, let us denote 
\begin{equation*}
  F(\vec{v};\beta)=:F_{\phi^\theta}^J(\alpha,\beta) . 
\end{equation*}
These functions are itself interesting objects. Actually,
for $n=2$, they define theta functions \cite{PoZa}.
We construct these explicitly in subsection \ref{ssec:m2}. 
The case $n=3$, they are related to Jacobi's indefinite theta functions, 
see \cite{pasol-poli:indefinite}. 
The $A_\infty$-relations then correspond to various identities between 
these generalized theta functions. 
A particular class of $A_\infty$-products and their properties are 
discussed in \cite{Poli:Eisen1, Poli:Eisen2}. 
As we see in Lemma \ref{lem:multi-beta}, 
the generalizations $F(\vec{v},{\frak B})$ are essentially obtained by
derivatives of $F_{\phi^\theta}^J(\alpha,\beta)$. 

These $F_{\phi^\theta}^J$ are $C^\infty$ at generic $(\alpha,\beta)$,
but they may not be continuous at some points
where CC-lines $\ti{v}$ arise. 
The terms $R^d$
we introduced in subsection \ref{ssec:str-const} 
correct the value at these points. 

Anyway, since 
these functions $F_{\phi^\theta}^J$ become
the main part of the structure constants of the $A_\infty$-structure,
they have various symmetries associated to
various $A_\infty$-automorphisms on $\cC$ or $\ti{\cC}$. 
In this section, we list some of such symmetries. 
We also discuss a $\theta$ independence of
$F_{\phi^\theta}^J$ in subsection \ref{ssec:theta-indep}, 
so finally we can drop $\theta$ and express them as $F_\phi^J$
for any $\theta$.

 \subsection{Translation invariance}
\label{ssec:translation-inv}
 
The Fukaya category $\cC_\theta$ is invariant under
the translation $(x,y)\mapsto (x,y+c)$ with $c\in\R$. 
Each object $a=(q_a,p_a,\alpha_a,\beta_a)$ is then
mapped to $(q_a,p_a,\alpha_a+c,\beta_a)$. This map is naturally lifted
to the linear $A_\infty$-automorphism $\cC\to\cC$,
where the structure constants of the $A_\infty$-structure are kept unchanged. 
This implies that 
\begin{equation*}
 F_{\phi^\theta}^J(\alpha+c,\beta)=  F_{\phi^\theta}^J(\alpha,\beta)
\end{equation*}
holds, where $\alpha+c:=(\alpha_{a_1}+c,\dots,\alpha_{n+1}+c)$. 

More generally, there exists the linear $A_\infty$-automorphism
corresponding to the translation
$(x,y)\mapsto (x+c_x,y+c_y)$, by which $q_a y= p_a x+\alpha_a$
is mapped to $q_a (y-c_y)= p_a (x-c_x) +\alpha_a$, i.e.,
$\alpha_a\mapsto \alpha_a+ (q_a c_y-p_a c_x)$.  
Thus, the function $F_{\phi^\theta}^J$ is constant along these two directions.

 \subsection{Periodicity}
\label{ssec:periodicity}
 
Any object $a=(q_a,p_a,\alpha_a,\beta_a)$ is isomorphic to
$a':=(q_a,p_a,\alpha_a+1,\beta_a)$. 
In particular, we see that 
the structure constant is invariant under this change $a\mapsto a'$. 
Thus, the function $F_\phi^J$ is periodic in $\alpha$, i.e.,
\begin{equation*}
  F_{\phi^\theta}^J(\alpha_{a_1},\dots\alpha_{a_i}+1,\dots,\alpha_{n+1},\beta)
  =  F_{\phi^\theta}^J(\alpha_{a_1},\dots\alpha_{a_i},\dots,\alpha_{n+1},\beta) 
\end{equation*}
for each $i=1,\dots,n+1$.

On the other hand, $F_{\phi^\theta}^J$ is not periodic in $\beta$. 
The behavior of $F_{\phi^\theta}^J$ under the change 
$(\beta_{a_1},\dots\beta_{a_i},\dots,\beta_{n+1})
\mapsto (\beta_{a_1},\dots\beta_{a_i}+1,\dots,\beta_{n+1})$ 
is also interesting from the viewpoint of automorphic functions. 
However, it depends on $n$ and $a_1,\dots,a_{n+1}$, so
we skip to examine the property here.

 \subsection{Cyclicity}
\label{ssec:cyclic-sym}

It is clear that $F_{\phi^\theta}^J$ has the cyclic symmetry
since this is the structure constant of a cyclic
$A_\infty$-product. 
By replacing $(1,\dots,n+1)$ by $(2,\dots,n+1,1)$,
all $\phi^\theta$, $J$, $\alpha$ and $\beta$ are rotated
simultaneously, and the resulting function $F$ is
the same as the original one or minus the original one.

 \subsection{$SL(2,\Z)$ invariance}
\label{ssec:SL2Z}
 
For an element $g\in SL(2,\Z)$, we consider the automorphism $g:T^2\to T^2$ 
defined by
\begin{equation*}
 \bp x\\ y\ep \mapsto  g\bp x\\ y\ep . 
\end{equation*} 
This preserves the constant symplectic form on $T^2$.
For each object $a=(q_a,p_a,\alpha_a,\beta_a)$,
we obtain the pullback $g^*(a)$.
Also, for $\theta$, we obtain $g^*(\theta)$. 
They are given explicitly by
\begin{equation*}
  \bp 1 & \theta \\
  -p_a & q_a \ep
  \bp x\\ y\ep
  \mapsto 
  \bp 1 & g^*(\theta) \\
  -g^*(p_a) & g^*(q_a) \ep
  \bp x\\ y\ep
 =   \bp 1 & \theta \\
  -p_a & q_a \ep
  g \bp x\\ y\ep
\end{equation*}
and $g^*(a)=(g^*(q_a),g^*(p_a),\alpha_a,\beta_a)$. 
This further
induces a linear $A_\infty$-automorphism $\cC_\theta\to\cC_{g^*(\theta)}$,
where the structure constants of the $A_\infty$-structure are kept unchanged. 
Thus, by denoting $g^*(\phi^\theta)$ the corresponding change of $\phi^\theta$, 
one has
\begin{equation*}
   F^J_{g^*(\phi^\theta)}=F^J_{\phi^\theta} . 
\end{equation*}

 \subsection{$\theta$ independence of $\ti{\cC}_\theta$}
\label{ssec:theta-indep}
 
Let us consider the additive $A_\infty$-category $\ti{\cC}_\theta$ of 
the Fukaya category $\cC_\theta$. 
We see that $\ti{\cC}_\theta$ is independent of $\theta$ 
\cite[Proposition 5.22]{hk:rims07nc} as follows. 
Any indecomposable object in $\ti{\cC}_\theta$ is
of the form $a[r_a]$, where $a\in\cC_\theta$ and $r_a\in\Z$. 
To such an object, we attach
\begin{equation*}
 \phi^\theta_{a[r_a]}:= \phi^\theta_a + r_a\pi . 
\end{equation*}
Take another $\theta'$ and the full subcategory of $\ti{\cC}_\theta$ 
consisting of all indecomposable objects $a[r_a]$ satisfying
\begin{equation*}
  \tan^{-1}\left(\frac{1}{\theta'}\right)-\pi< \phi^\theta_{a[r_a]}
  \le \tan^{-1}\left(\frac{1}{\theta'}\right)
\end{equation*}
is exactly $\cC_{\theta'}$ itself. 
Of course, each object $a[*]\in\cC_{\theta'}$ is included as
$a[**]\in\ti{\cC}_\theta$ 
so that $\phi^{\theta'}_{a[*]}=\phi^\theta_{a[**]}$. 
This fact corresponds to the mirror dual statement of what is discussed
in \cite{PoSc}. 

Thus, we denote $\ti{\cC}_\theta=:\ti{\cC}$. 
In $\ti{\cC}$, we can define
a CC-collection of indecomposable objects $(a_1[r_1],\dots,a_{n+1}[r_{n+1}])$ 
which satisfies
\begin{equation*}
   t-\pi < \phi_{a_i[r_i]} \le t
\end{equation*}
for all $i\in\{1,\dots,n+1\}$ with some $t\in\R$.
Then, we can define the corresponding functions 
$F_{\phi}^J$ with $\phi=(\phi_1,\dots,\phi_{n+1})$ so that these satisfy
\begin{equation*}
  F_{\phi'}^J=(-1)^{nr}F_{\phi}^J
\end{equation*}
for 
\begin{equation*}
  \phi'=(\phi'_{a_1},\dots,\phi_{a_{n+1}}')
  =(\phi_{a_1}+r\pi,\dots,\phi_{a_{n+1}}+r\pi) . 
\end{equation*}
Namely, we can drop $\theta$ in $\phi^\theta$ 
even if we start with $\cC_\theta$ with some fixed $\theta$. 
Note that the sign $(-1)^{nr}$ comes from
\begin{equation}\label{str-const-sign}
  \ti\eta(\ti{m}_n(w'_{12},\dots,w'_{n(n+1)}),w'_{(n+1)1})
  =(-1)^{(|\eta|+1)r_1+nr_1}
  \ti\eta(\ti{m}_n(w_{12},\dots,w_{n(n+1)}),w_{(n+1)1})
\end{equation}
for $w_{ij}\in\ti\cC(a_i,a_j)$ and $w'_{ij}\in\ti\cC(a_i[r_i],a_j[r_j])$, 
which is obtained by our settings 
(\ref{tim}) and (\ref{tieta}) with $|\eta|=1$.

For $t<t'$ such that $(t-\pi,t]\cap (t'-\pi,t']\neq \emptyset$, 
let us consider $F_\phi^J$ such that 
$t-\pi<\phi_i\le t$ for all $i=1,\dots,n+1$. 
For each $i$, we set
\begin{equation*}
  \phi'_i:=
  \begin{cases}
    \phi_i+\pi & t-\pi < \phi_i \le t'-\pi \\
    \phi_i &  t'-\pi< \phi_i\le t 
  \end{cases} . 
\end{equation*}
Then, $t'-\pi< \phi'_i\le t'$ holds for all $i$. 
In this situation, one has
\begin{equation*}
  F_{\phi'}^J = (-1)^{n(\phi'_1-\phi_1)/\pi} F_\phi^J
\end{equation*}
again by (\ref{str-const-sign}).

\section{Examples}
\label{sec:exmp}

In this section, we discuss a few examples of the structure constants of
the $A_\infty$-products in $\cC$. 
In subsection \ref{ssec:m2}, we calculate the structure constants of
all transversal products $m_2$. 
In subsection \ref{ssec:m3},
we calculate a few non-transversal triple products and discuss
their relation to an exact triangle in $\Tri(\cC)$.

 \subsection{The products $m_2$}
\label{ssec:m2}

In this subsection, we give the formula for 
the structure constants $\varphi(v_{ab}^{j_{ab}},v_{bc}^{j_{bc}},v_{ca}^{j_{ca}})$ 
for any objects $a,b,c$. 

As a practice, we first start with doing it for the special case
where 
\begin{equation*}
  \begin{split}
 & c=(q_c,p_c,\alpha_c,\beta_c)=(1,1,0,0) \\
 & a=(q_a,p_a,\alpha_a,\beta_a)=(1,0,0,0) \\
 & b=(q_b,p_b,\alpha_b,\beta_b)=(2,1,\alpha,\beta)  
  \end{split}
\end{equation*}
with fixed $\alpha$ and $\beta$. 
This situation is related to the example 
discussed in the next subsection.

In this case, all the vector space $V_{ab}$, $V_{bc}$,
$V_{ca}$ are one dimensional fortunately. 
We denote the corresponding bases by $v_{ab}$, $v_{bc}$ and $v_{ca}$,
where $|v_{ab}|=|v_{bc}|=0$ and $|v_{ca}|=1$. 

We fix $v_{ca}$ in $\R^2$, which is the intersection point of
$L_c$ and $L_a$, where 
\begin{equation*}
  \begin{split}
    L_c: y= x, \\
    L_a: y= 0 .
  \end{split}
\end{equation*}
We further consider the line
\begin{equation*}
  \ti{L}_b: 2 y =  x + \alpha + m
\end{equation*}
in $\R^2$ with $m\in\Z$, and define
$\ti{v}_{ab}$ and $\ti{v}_{bc}$ as the intersection points of
$(L_a,\ti{L}_b)$ and $(\ti{L}_b,L_c)$, respectively. 
One has 
\begin{equation*}
  (x,y)(\ti{v}_{ab})= (-\alpha-m,0),
  \qquad (x,y)(\ti{v}_{bc})=(\alpha+m,\alpha+m) . 
\end{equation*}
The area of the triangle
$\ti{v}:=\ti{v}_{ab}\ti{v}_{bc}v_{ca}$ is then $(\alpha+m)^2$,
and $\Delta_b(\ti{v})=\alpha+m$. 
Thus, we obtain 
\begin{equation}\label{m2}
 \varphi(v_{ab},v_{bc},v_{ca}) = 
 \sum_{m\in\Z} \exp(2\pi\ii\rho (\alpha+m)^2 + 2\pi\ii\beta(\alpha+m))  . 
\end{equation}

For general $a,b,c$, the situation is more complicated
because the dimension of $V_{ab}$, $V_{bc}$, $V_{ca}$ may not
be one. 
We employ the arguments in subsection \ref{ssec:basis} fully, 
and consider the following lines in $\wt{T^2}$: 
\begin{align}
 L_c:\ \ &q_c y=p_c x+\alpha_c\ ,\\
 L'_a:\ \ &q_a y=p_a x-q_{ca}j_{ca}+\alpha_a\ ,\\  
 \ti{L}'_b:\ \ &q_b y=p_b x-q_{cb}j_{ca}
 -q_{ab}j_{ab}-p_{ab}m+\alpha_b\ .
\end{align}
Here, the intersecting point of $L_c$ with $L'_a$ is $v_{ca}$, 
and $\ti{L}'_b$ is defined so that the intersecting point
$\ti{v}_{ab}^{j_{ab}}$ of $L'_a$ with $L'_b$ satisfies
$\pi(\ti{v}_{ab}^{j_{ab}})=\pi(v_{ab}^{j_{ab}})$. 
These facts in particular mean that 
$\pi(L_a')=\pi(L_a)$ and $\pi(L_b')=\pi(L_b)$ hold. 

More intuitively, $L'_a$ is defined so that it passes through 
$$
(x,y)=
\left(q_c\frac{q_{ca}j_{ca}-\alpha_a}{p_{ca}},
p_c\frac{q_{ca}j_{ca}-\alpha_a}{p_{ca}}\right)\ ,
$$ 
and $L'_b$ is defined so that it passes through 
$$
(x,y)=\left(
q_c\frac{q_{cb}j_{ca}}{p_{cb}}+q_a\left(m+\frac{q_{ab}j_{ab}-\alpha_b}{p_{ab}}
\right), 
p_c\frac{q_{cb}j_{ca}}{p_{cb}}+p_a\left(m+\frac{q_{ab}j_{ab}-\alpha_b}{p_{ab}}
\right)
\right)\ .
$$ 

Now, one has 
\begin{align}
 &(x, y)(\ti{v}_{ca}^{j_{ca}})
=\left(\frac{q_cq_{ca}j_{ca}-q_c\alpha_a+q_a\alpha_c}{p_{ca}}, 
\frac{p_cq_{ca}j_{ca}-p_c\alpha_a+p_a\alpha_c}{p_{ca}}\right) 
\label{xy-eab}\\
 &(x,y)(\ti{v}_{ab}^{j_{ab}})
 =\left(\frac{q_ap_{ab}m
+(q_aq_{cb}-q_bq_{ca})j_{ca}+q_aq_{ab}j_{ab}
-q_a\alpha_b+q_b\alpha_a}{p_{ab}}, \right.\\
&\hspace*{4.0cm} \left.
\frac{q_ap_{ab}m+(p_aq_{cb}-p_bq_{ca})j_{ca}+p_aq_{ab}j_{ab}
-p_a\alpha_b+p_b\alpha_a}{p_{ab}}\right)\ .
 \label{xy-ebc}
\end{align}
Since the slopes of the three lines are fixed, 
the triangle is determined when the two vertices $v_{ca}^{j_{ca}}$ and 
$\ti{v}_{ab}^{j_{ab}}$ are given. 
Namely, the triangle exists 
only when the intersection point $\ti{v}_{bc}$ of $(\ti{L}'_b,L_c)$ 
satisfies $\pi(\ti{v}_{bc})=\pi(v_{bc}^{j_{bc}})$.
Since 
\begin{equation*}
 (x,y)(\ti{v}_{bc})=\left(
\frac{q_c(p_{ab}m+q_{cb}j_{ca}+q_{ab}j_{ab})-q_c\alpha_b+q_b\alpha_c}
{p_{cb}}, 
\frac{p_cp_{ab}m+p_cq_{cb}j_{ca}+p_cq_{ab}j_{ab}-p_c\alpha_b+p_b\alpha_c}
{p_{cb}}\right)\ 
\end{equation*}
and 
\begin{equation*}
 (x,y)(\ti{v}_{cb}^{j_{cb}})
 =\left(
\frac{q_cq_{cb}j_{cb}-q_c\alpha_b+q_b\alpha_c}{p_{cb}}, 
\frac{p_bq_{cb}j_{cb}-p_c\alpha_b+p_b\alpha_c}{p_{cb}}
\right)\ ,
\quad -q_{cb}j_{cb}=j_{bc} , 
\end{equation*}
the condition $\pi(\ti{v}_{bc})=\pi(v_{bc}^{j_{bc}})$ is equivalent to 
\begin{equation*}
  \begin{split}
  & q_c(p_{ab}m+q_{cb}j_{ca}+q_{ab}j_{ab}+j_{bc})\in p_{bc}\Z, \\
  & p_c(p_{ab}m+q_{cb}j_{ca}+q_{ab}j_{ab}+j_{bc})\in p_{bc}\Z, 
  \end{split}
\end{equation*}
but we see that the above two identities are equivalent to the single
identity 
\begin{equation*}
 p_{ab}m+q_{cb}j_{ca}+q_{ab}j_{ab}+j_{bc} \in p_{bc}\Z\ 
\end{equation*}
since $q_c$ and $p_c$ are relatively prime. 
Let us define the corresponding Kronecker's delta:
\begin{equation*}
  \delm{p_{bc}}_i^j:=
  \begin{cases}
    1 & i-j\in p_{bc}\Z \\
    0 & \text{otherwise}\ .
  \end{cases}
\end{equation*}
The structure constant is then given by 
\begin{equation}
  \varphi(v_{ab}^{j_{ab}},v_{bc}^{j_{bc}},v_{ca}^{j_{ca}})
 =\sum_{m\in\Z}\delm{p_{bc}}_{p_{ab}m+q_{cb}j_{ca}+q_{ab}j_{ab}+j_{bc}}^0
 \exp{\left(2\pi\ii\rho\triangle_m\right)}\exp{(2\pi\ii \int\beta_m)}\ ,
 \label{strconst2}
\end{equation}
where $\triangle_m$ is the area of the triangle made from 
$L_c$, $L'_a$ and $\ti{L}'_b$ which is computed as 
\begin{equation}
 \begin{split}
 & \triangle_m
=\ov{-2p_{ab}p_{bc}p_{ca}}
\left(l_{abc}(m,j_{ab},j_{bc},j_{ca})-\alpha_{abc}\right)^2\ ,\\
& l_{abc}(m,j_{ab},j_{bc},j_{ca}):=p_{ca}p_{ab}m+p_{ca}q_{ab}j_{ab}-p_{ab}j_{ca}\ ,\\
& \alpha_{abc}:=p_{bc}\alpha_a+p_{ca}\alpha_b+p_{ab}\alpha_c , 
 \end{split}
 \label{tri}
\end{equation}
and 
\begin{equation*}
 \begin{split}
  & \exp{(2\pi\ii \int\beta_m)}:=\exp{\left(-2\pi\ii
 \beta_{abc}
\left(l_{abc}(m,j_{ab},j_{bc},j_{ca})-\alpha_{abc}\right)\right)}\ , \\
  & \beta_{abc}:= p_{bc}\beta_a+p_{ca}\beta_b+p_{ab}\beta_c . 
 \end{split}
\end{equation*}
Though the cyclicity of
this structure constant 
$\varphi(v_{ab}^{j_{ab}},v_{bc}^{j_{bc}},v_{ca}^{j_{ca}})$ is guaranteed by
the homological mirror symmetry, it can be checked directly
in the explicit formula above; since $\alpha_{abc}$ and $\beta_{abc}$
are already cyclic, we may check the cyclicity of $l_{abc}$ under the
constraint coming from $\delm{p_{bc}}$, see \cite{hk:nchms}.

 \subsection{Certain triple products and an exact triangle}
\label{ssec:m3}

The triangulated category $\Tri(\cC)$ includes 
informations of $A_\infty$-products in $\cC$.  
In this subsection, we discuss 
a relation between $A_\infty$-products of the Fukaya category $\cC$
and an exact triangle in $\Tri(\cC)$ by an example. 
The exact triangle we discuss below is the mirror dual of
that discussed in \cite{kobayashi:2tori}. 
The reader can confirm that various complicated calculations in
\cite{kobayashi:2tori} are already 
included in obtaining the Fukaya category $\cC$. 

The objects we treat are those considered in the previous subsection. 
For 
\begin{equation*}
 \begin{split}
 & a=(q_a,p_a,\alpha_a,\beta_a)=(1,0,0,0), \\
 & c=(q_c,p_c,\alpha_c,\beta_c)=(1,1,0,0), 
 \end{split}
\end{equation*} 
$V_{ca}$ is a one-dimensional vector space spanned by
degree one base $v_{ca}$.
Regarding it as a map $v_{ca}\in\ti{\cC}(c[-1],a)$, we consider
the cone 
\begin{equation*}
  C(v_{ca}):=(c\oplus a, \left(\bps 0 & \bar{v}_{ca} \\ 0 & 0\eps\right))
  \in \Tw(\cC) . 
\end{equation*}
Now, we will show that this one-sided twisted complex is
isomorphic in $\Tri(\cC)$ to $(b,0)$ with 
$b=(q_b,p_b,\alpha_b,\beta_b)=(2,1,1/2,1/2)$. 

Here, we first assume that $b=(2,1,\alpha_b,\beta_b)$ with
arbitrary $(\alpha_b,\beta_b)$.
The graded vector space
$\Tw(\cC)((b,0),C(v_{ca}))$ is 
spanned by the degree zero element 
\begin{equation}\label{g1}
  \bp v_{bc} & 0 \ep . 
\end{equation}
Also, $\Tw(\cC)(C(v_{ca}),(b,0))$ is 
spanned by the degree zero element 
\begin{equation}\label{g2}
  \bp 0 \\ v_{ab} \ep . 
\end{equation}
The space $\Tri(\cC)((b,0),C(v_{ca}))$ is the zero-th
cohomology of the complex
$(\Tw(\cC)((b,0),C(v_{ca})),m_1^\Tw)$.
The degree-zero element $\left(\bps v_{bc} & 0 \eps\right)$ 
is then closed, 
\begin{equation*}
  m_1^\Tw\left( \bp \bar{v}_{bc} & 0 \ep \right)) =
  \ti{m}_2\left(\bp 0 \ep , \bp \bar{v}_{bc} & 0\ep\right)
  + \ti{m}_2\left(\bp \bar{v}_{bc} & 0\ep , 
 \bp 0 & \bar{v}_{ca} \\ 0 & 0\ep \right) =0 , 
\end{equation*}
if and only if the product $m_2(v_{bc},v_{ca})$ is zero.
Similarly, 
$\left(\bps 0 \\ v_{ab} \eps\right)$ is closed,
\begin{equation*}
  m_1^\Tw\left(\bp 0 \\ \bar{v}_{ab} \ep\right)  = 
  \ti{m}_2\left(
  \bp 0 & \bar{v}_{ca} \\ 0 & 0\ep , 
  \bp 0 \\ \bar{v}_{ab} \ep 
   \right)
 + \ti{m}_2\left(
  \bp 0 \\ \bar{v}_{ab} \ep , 
 \bp 0 \ep
  \right) =0 , 
\end{equation*}
if and only if $m_2(v_{ca},v_{ab})=0$. 

One can understand these facts without direct calculations as above. 
Applying the cohomological functor $\Tri(\cC)((b,0),\ \ )$
to the exact triangle 
\begin{equation}\label{et}
 \cdots \to (c[-1],0)\to (a,0) \to C(v_{ca})\to (c,0)\overset{v_{ca}}{\to} (a[1],0)\to C(v_{ca})[1]\to (c[1],0)\to\cdots 
\end{equation}
yields the long exact sequence
\begin{equation*}
  \cdots\to 0 \to 0\to \Tri(\cC)((b,0),C(v_{ca}))\to \la v_{bc}\ra
  \overset{v_{ca}}{\to}
  \la v_{ba}\ra\to\Tri(\cC)((b,0),C(v_{ca})[1])\to 0\to\cdots \ .
\end{equation*}
Here, 
if the map $v_{ca}=\Tri(\cC)((b,0),v_{ca})$ in the exact sequence above
is non-zero, 
then $\Tri(\cC)((b,0),C(v_{ca}))$ turns out to be trivial and
so $(b,0)$ cannot be isomorphic to $C(v_{ca})$. 
Now, note that the image of $v_{bc}$ by the map $v_{ca}=\Tri((b,0),v_{ca})$
is $m_2(v_{bc},v_{ca})$. 
Thus, one sees $\Tri(\cC)((b,0),C(v_{ca}))=0$ if $m_2(v_{bc},v_{ca})\ne 0$,
and $\Tri(\cC)((b,0),C(v_{ca}))\simeq\C$ if $m_2(v_{bc},v_{ca})=0$.
Parallel fact is obtained for $\Tri(\cC)(C(v_{ca}),(b,0))$ 
by applying the cohomological functor $\Tri(\cC)(\ \ , (b,0))$
to the exact triangle (\ref{et}). 

Next, assume that $m_2(v_{bc},v_{ca})=0$
and then check when 
the generators (\ref{g1}) and (\ref{g2}) form isomorphisms.
One has
\begin{equation*}
 \begin{split}
  m_2^\Tw\left(\bp \bar{v}_{bc} & 0 \ep , \bp 0 \\ \bar{v}_{ab} \ep\right)
  & =\sum_i\ti{m}_{2+i} \left(\bp \bar{v}_{bc} & 0 \ep ,
  \bp 0 & \bar{v}_{ca} \\ 0 & 0\ep^i ,\bp 0 \\ \bar{v}_{ab} \ep  \right) \\
  & = \ti{m}_3\left(\bp \bar{v}_{bc} & 0 \ep ,
  \bp 0 & \bar{v}_{ca} \\ 0 & 0\ep ,\bp 0 \\ \bar{v}_{ab} \ep  \right) \\
  & = \bar{m}_3(\bar{v}_{bc},\bar{v}_{ca},\bar{v}_{ab}) . 
 \end{split}
\end{equation*}
and 
\begin{equation*}
 \begin{split}
m_2^\Tw\left(\bp 0 \\ \bar{v}_{ab} \ep ,    \bp \bar{v}_{bc} & 0 \ep\right)
& =\sum_{i,j}\ti{m}_{2+i+j} \left(
\bp 0 & \bar{v}_{ca} \\ 0 & 0\ep^i , \bp 0 \\ \bar{v}_{ab}\ep,
\bp \bar{v}_{bc} & 0 \ep , 
  \bp 0 & \bar{v}_{ca} \\ 0 & 0\ep^j  \right) \\
  & = 
  \bp
  \bar{m}_3(\bar{v}_{ca},\bar{v}_{ab},\bar{v}_{bc}) & 0 \\
  \bar{m}_2(\bar{v}_{ab},\bar{v}_{bc}) &
  \bar{m}_3(\bar{v}_{ab},\bar{v}_{bc},\bar{v}_{ca})
  \ep . 
 \end{split}
\end{equation*}
Since $m_2(v_{ab},v_{bc})=0$ by the cyclicity
$0=\omega(m_2(v_{bc},v_{ca}),v_{ab})=-\omega(m_2(v_{ab},v_{bc}),v_{ca})$,
the generators (\ref{g1}) and (\ref{g2}) form isomorphisms
if and only if
the structure constants $c_1,c_2,c_3$ of 
$\bar{m}_3(\bar{v}_{bc},\bar{v}_{ca},\bar{v}_{ab})=c_1\1_b$,
$\bar{m}_3(\bar{v}_{ca},\bar{v}_{ab},\bar{v}_{bc})=c_2\1_c$,
$\bar{m}_3(\bar{v}_{ab},\bar{v}_{bc},\bar{v}_{ca})=c_3\1_a$
satisfy $c_1=c_2=c_3\ne 0$. 

Now, let us discuss when these conditions are satisfied. 
First, the structure constant
$\varphi(v_{ab},v_{bc},v_{ca})$ is already calculated in
(\ref{m2}).
Let us express this as 
\begin{equation*}
 \begin{split}
  & \varphi(v_{ab},v_{bc},v_{ca})= \sum_{m\in\Z} F_m(\alpha,\beta) \\
   & \ \ F_m((\alpha,\beta)):=
   \exp(2\pi i\rho (\alpha+m)^2 + 2\pi i\beta (\alpha+m)) . 
 \end{split}
\end{equation*}
When $\alpha=\beta=1/2$, we see that 
$F_m(1/2,1/2)= - F_{-m-1}(1/2,1/2)$ and 
we have $\varphi(v_{ab},v_{bc},v_{ca})= 0$. 
We can also see this fact because
$\sum_{m\in\Z} F_m(\alpha,\beta)$ is a theta function and
$(\alpha,\beta)=(1/2,1/2)$ is a zero of it. 
Thus, for $\alpha=\beta=1/2$, we have
\begin{equation*}
   m_2(v_{ab},v_{bc})=0,\quad m_2(v_{bc},v_{ca})=0,\quad m_2(v_{ca},v_{ab})=0 . 
\end{equation*}

Next, let us consider
the triple product
$m_3(v_{ab},v_{bc},v_{ca})$. 
By Lemma \ref{lem:multi-beta}, we have 
\begin{equation*}
   m_3(v_{ab},v_{bc},v_{ca})
   = \frac{\sigma(a)}{2\pi\ii}\left(\fpartial{\beta}\sum_{m\in\Z}  F_m(1/2,\beta) |_{\beta=1/2}\right)\cdot 1_a, \qquad \sigma(a)=1  , 
\end{equation*}
in $\cC$. 
Rewriting this as an $A_\infty$ triple product in $s\cC$
by (\ref{m-suspension}), 
it turns out that 
\begin{equation*}
  c_1=-\frac{1}{2\pi\ii}\fpartial{\beta}
  \left(\sum_{m\in\Z}  F_m(1/2,\beta) |_{\beta=1/2}\right) . 
\end{equation*}
One can similarly calculate $c_2$ and $c_3$, and obtain 
$c_1=c_2=c_3$.

Since $\sum_{m\in\Z} F_m(\alpha,\beta)$ is the theta function,
the zeros of its derivative is known 
by the Jacobi's derivative formula. 
Actually $\alpha=\beta=1/2$ is not a zero, 
so one has $c_1=c_2=c_3\ne 0$. See \cite{kobayashi:2tori}. 

As a result, one has $C(v_{ca})\simeq b=(2,1,1/2,1/2)$.
Since $b'=(2,1,\alpha,\beta)\simeq b$ if and only if
$\alpha-1/2\in\Z$ and $\beta-1/2\in\Z$ hold, we can conclude that 
\begin{equation*}
C(v_{ca})\simeq b=(2,1,\alpha,\beta)
\end{equation*}
if and only if $\alpha-1/2\in\Z$ and $\beta-1/2\in\Z$ hold.

\end{document}